\documentclass[11pt]{article}
\usepackage[utf8]{inputenc}
\usepackage{fullpage,amsmath,amssymb,bm,url,mathrsfs}
\usepackage{algorithmic}
\usepackage{algorithm}
\usepackage{amsthm}
\usepackage{hyperref}
\usepackage{makecell}
\usepackage{tablefootnote}
\usepackage{comment}
\usepackage{soul}
\usepackage{graphics,graphicx}
\usepackage{subcaption,caption,cleveref}
\hypersetup{
	colorlinks=true,
	linkcolor=blue,
	citecolor=blue,
	filecolor=magenta,      
	urlcolor=cyan,
}

\usepackage[round,colon,authoryear]{natbib}
\usepackage{xcolor}
\usepackage{dsfont}
\usepackage{bbm}
\usepackage{enumitem}
\usepackage{wrapfig}
\usepackage{relsize}

\def\calA{{\mathcal A}}

\def\calF{{\mathcal F}}

\def\calL{{\mathcal L}}
\def\calM{{\mathcal M}}

\def\calP{{\mathcal P}}

\def\calR{{\mathcal R}}
\def\calS{{\mathcal S}}

\def\calU{{\mathcal U}}

\def\calX{{\mathcal X}}

\def\BB{{\mathbb B}}

\def\EE{{\mathbb E}}

\def\II{{\mathbb I}}

\def\OO{{\mathbb O}}

\def\RR{{\mathbb R}}
\def\SS{{\mathbb S}}

\def\bmu{{\boldsymbol \mu}}

\DeclareMathOperator*{\argmin}{arg\,min}

\def\calA{{\cal  A}} \def\scrA{{\mathscr  A}}
 
 \def\scrC{{\mathscr  C}}
 
 \def\scrE{{\mathscr  E}}
\def\calF{{\cal  F}} \def\scrF{{\mathscr  F}}
 \def\scrG{{\mathscr  G}}
 
 \def\scrI{{\mathscr  I}}
 
 \def\scrK{{\mathscr  K}}
\def\calL{{\cal  L}} 
\def\calM{{\cal  M}} \def\scrM{{\mathscr  M}}

\def\calP{{\cal  P}} \def\scrP{{\mathscr  P}}
 
\def\calR{{\cal  R}} 
\def\calS{{\cal  S}} \def\scrS{{\mathscr  S}}
 \def\scrT{{\mathscr  T}}
\def\calU{{\cal  U}} \def\scrU{{\mathscr  U}}

\def\calX{{\cal  X}} \def\scrX{{\mathscr  X}}
 \def\scrY{{\mathscr  Y}}
 \def\scrZ{{\mathscr Z}}

\newcommand{\bfm}[1]{\ensuremath{\mathbf{#1}}}

\def\ba{\bfm a}   \def\bA{\bfm A}  
\def\bb{\bfm b}   \def\bB{\bfm B}  \def\BB{\mathbb{B}}
\def\bc{\bfm c}   \def\bC{\bfm C}  
     
\def\be{\bfm e}   \def\bE{\bfm E}  \def\EE{\mathbb{E}}

   \def\bI{\bfm I}  \def\II{\mathbb{I}}

   \def\bM{\bfm M}  
     
     \def\OO{\mathbb{O}}
\def\bp{\bfm p}     
     
     \def\RR{\mathbb{R}}
   \def\bS{\bfm S}  \def\SS{\mathbb{S}}
     
\def\bu{\bfm u}   \def\bU{\bfm U}  
\def\bv{\bfm v}   \def\bV{\bfm V}  
\def\bw{\bfm w}     
\def\bx{\bfm x}   \def\bX{\bfm X}  
\def\by{\bfm y}     
\def\bz{\bfm z}   \def\bZ{\bfm Z}

\def\hat{\widehat}

\theoremstyle{remark}

\newcommand{\eps}{\varepsilon}

\def\eps{\varepsilon}

\newcommand{\vertiii}[1]{{\left\vert\kern-0.25ex\left\vert\kern-0.25ex\left\vert #1 
		\right\vert\kern-0.25ex\right\vert\kern-0.25ex\right\vert}}

\usepackage{tikz}
\usetikzlibrary{decorations.pathreplacing,calc, patterns}

\def\scrE{\mathscr{E}}
\def\scrX{\mathscr{X}}
\def\scrT{\mathscr{T}}

\def\tilde{\widetilde}
\def\hat{\widehat}

\begin{document}

	\title{Tensor Methods in High Dimensional Data Analysis: Opportunities and Challenges}
	
	\author{Arnab Auddy$^\ast$, Dong Xia$^\dagger$ and Ming Yuan$^\ddagger$\\
	$^\ast$University of Pennsylvania\\
	$^\dagger$Hong Kong University of Science and Technology\\
	$^\ddagger$Columbia University
	}

	\date{(\today)}
	
\footnotetext[1]{
	Arnab Auddy's research was supported by NSF Grant DMS-2015285.}
\footnotetext[2]{
	Dong Xia's research was supported by Hong Kong RGC grant GRF 16302020 and GRF 16300121.}
\footnotetext[3]{
	Ming Yuan's research was supported by NSF Grants DMS-2015285 and DMS-2052955.}
\maketitle
	
\begin{abstract}
	Large amount of multidimensional data represented by multiway arrays or tensors are prevalent in modern applications across various fields such as chemometrics, genomics, physics, psychology, and signal processing. The structural complexity of such data provides vast new opportunities for modeling and analysis, but efficiently extracting information content from them, both statistically and computationally, presents unique and fundamental challenges. Addressing these challenges requires an interdisciplinary approach that brings together tools and insights from statistics, optimization and numerical linear algebra among other fields. Despite these hurdles, significant progress has been made in the last decade. This review seeks to examine some of the key advancements and identify common threads among them, under eight different statistical settings. 
\end{abstract}


	\section{Introduction}\label{sec:intro}
	There is a long and rich history of the use of tensor methods in data analysis. See, e.g, \cite{coppi1989multiway, kroonenberg2008applied, mccullagh2018tensor} and references therein. Indeed, many of the most popular tensor decomposition methods originated from factor analysis of multiway data \citep[see, e.g.,][]{tucker1966some, carroll1970analysis, harshman1970foundations}. The use of tensor methods, however, have enjoyed a renaissance of sorts in recent years. See, e.g., \cite{cichocki2015tensor,sidiropoulos2017tensor,bi2021tensors}. On the one hand, this is the direct consequence of collecting more data with more complex relational structure. For example, it is much more intuitive to think of a video as three dimensional array than a matrix, and analyze it as such. On the other hand, methods of moments have emerged as the most popular approach to various high dimensional latent factor models for which the maximum likelihood estimates are often difficult to compute; and tensors arise naturally when considering higher order moments of a random vector.
	
	Following the convention in the related literature, here and in what follows, we shall use the term ``tensor'' loosely in that (i) we shall not make an effort to differentiate between tensors and multiway arrays and (ii) when referring to tensors, we assume they are of order at least three so that they can be differentiated from matrices or vectors, e.g., tensors of order two or one.
	
	\subsection{Why Tensor Methods}
	Tensor methods have presented numerous new opportunities for statistical modeling and data analysis that are not available for matrix-based techniques. Conceptually, they offer a better way to organize data with complex relational structures and lead to analysis with enhanced interpretability. From a modeling perspective, they ensure the identifiability of latent factors and other parameters that may not otherwise be well defined. Finally, it is also worthwhile to use tensor methods from an inferential point of view because the multilinear structure can help restrict the impact of sampling error. We now use principal component analysis (PCA) as an example to illustrate these key merits.
	
	\subsubsection{Conceptual Benefit -- Interpretability} Most of the statistical techniques are designed for data matrices consisting of observations (rows) for multiple variables (columns). They can be inadequate for analyzing data with richer and more complex structures, and addressing research questions such as ``Who Does What to Whom and When'' \citep{kroonenberg2008applied}.
	
	\begin{wrapfigure}[11]{r}{0.4\textwidth}
		\begin{center}
			\vskip -25pt
			\includegraphics[width=0.4\textwidth]{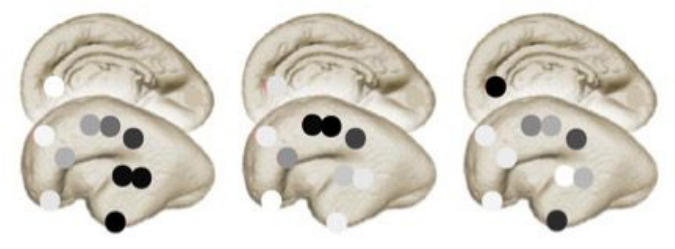}
		\end{center}
		\caption{Loadings of the first principal components for different brain regions: the darkness of the dots represent the magnitude of the loadings. Reproduced from \cite{liu2022characterizing}.}
		\label{fig:gene}
	\end{wrapfigure}
	
	Consider, for example, applying PCA to the spatiotemporal gene expression data of the human brain from \cite{kang2011spatio}. The data, after preprocessing, contains expression of 1087 genes over 10 different brain regions for 13 developmental periods. To apply the classical PCA, one would be forced to discard the spatial temporal information and format the data into a matrix of 1087 genes (rows) and 130 combinations of location and period (columns). On the other hand, as demonstrated by \cite{liu2022characterizing}, treating the data as third order tensor enables us to delineate the growth pattern and spatial distribution of gene expression. In particular, as shown in Figure \ref{fig:gene}, the factor loadings of the first principal components for different brain regions highlights the relevance of spatial proximity.
	
	\subsubsection{Benefit for Modeling -- Identifiability} A remarkable property of tensors, first discovered by \cite{kruskal1977three}, is that nearly all low rank tensors can be uniquely written as a sum of rank-one tensors. This is to be contrasted with the fact that there are infinite many ways to write a low rank matrix as a sum of rank-one matrices. The implication of such identifiability can be demonstrated in the comparison between PCA and independent component analysis (ICA). The principal components (PC) can be identified as the singular vectors of the covariance matrix, and likewise, the independent components (IC) are associated with the ``singular vectors'' of the fourth order cumulant tensor. Suppose that we observe $d$ independent random variables with mean zero, unit variance and nonzero excess kurtosis, up to an unknown rotation. We cannot reconstruct these variables via PCA because the variance is invariant to any rotation. We can however do so from its fourth order cumulant tensor which has rank $d$ with each component corresponding to one of the independent variables.
	
	\begin{figure}[htbp]
		\centering
		\begin{subfigure}[t]{0.5\textwidth}
			\centering
			\includegraphics[height=1.2in]{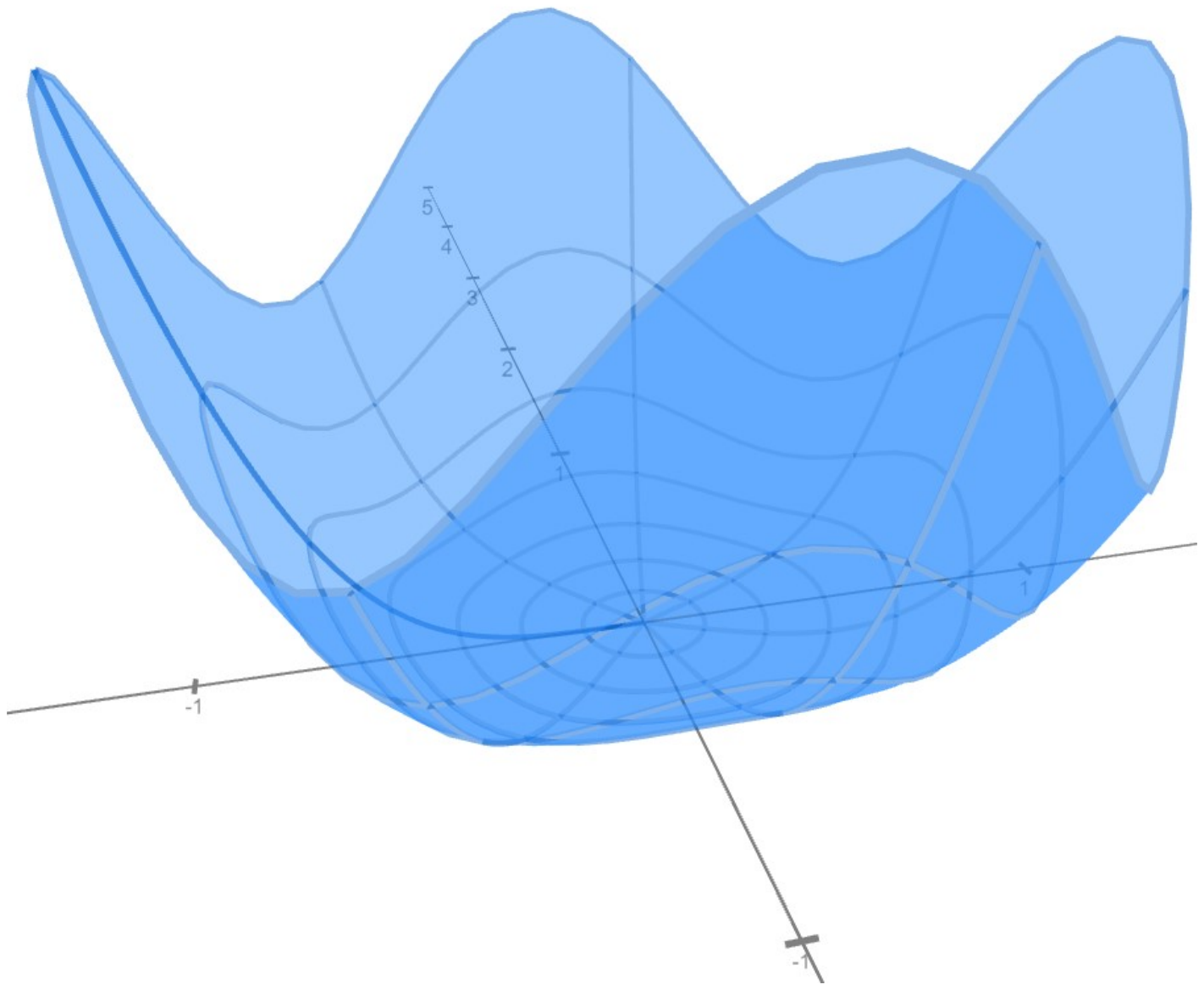}
		\end{subfigure}%
		~ 
		\begin{subfigure}[t]{0.5\textwidth}
			\centering
			\includegraphics[height=1.2in]{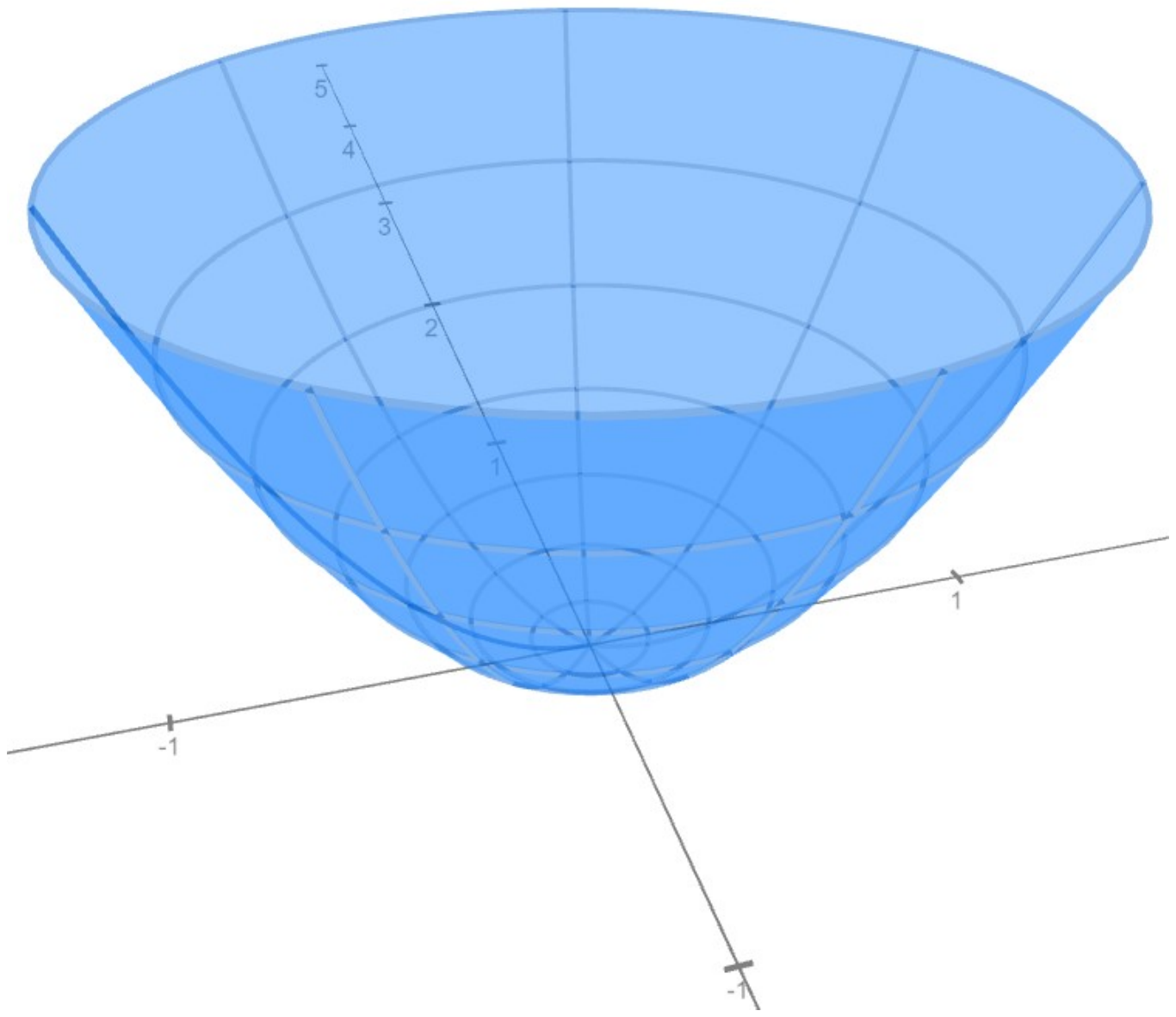}
		\end{subfigure}
		\caption{ICA vs PCA: $\kappa_4(\bu^\top \bX)=\EE(\bu^\top \bX)^4-3$ and $\texttt{cov}(\bu^\top \bX)$ as functions of $\bu$ over the unit circle.}
		\label{fig:ica}
	\end{figure}
	
	As a more concrete example, consider the case of two independent random variables $\bX=(X_1,X_2)^\top$ with $\texttt{var}(\bX)=\II_d$. The left panel of Figure \ref{fig:ica} shows the excess kurtosis as a function of $\bu=(u_1,u_2)^\top$: $f(\bu):=\kappa_4(u_1 X_1+u_2X_2)$ over the unit circle. There are only four maximizers corresponding to the independent components $\pm X_1$ and $\pm X_2$ respectively. In contrast, $\texttt{cov}(\bu^\top \bX)$ remains constant over all $\bu$ on the unit circle, as shown in the right panel. 
	
	\subsubsection{Benefit for Inferences -- Restricted Effect of Perturbation} Statistical analysis is based on noisy observations and subject to sampling error. Understanding the effect of sampling error plays a central role in any inference procedure. A recent observation is that, due to the inherent structural constraints, effect of sampling error on tensors can be more isolated than matrices which may allow us to conduct statistical inferences that are infeasible or impossible for data matrices.
	
	\begin{figure}[htbp]
		\centering
		\includegraphics[width=0.6\textwidth]{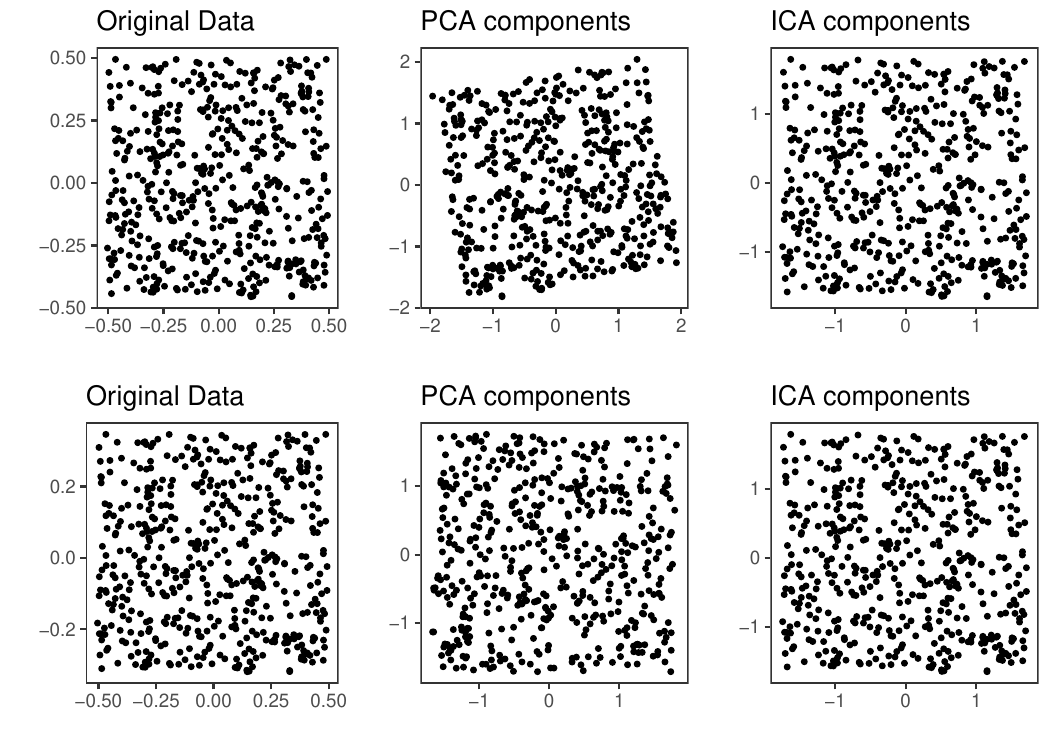}
		\caption{Effect of eigengap on PCA and ICA.}
		\label{fig:ica2}
	\end{figure}
	
	Take the classical PCA as an example. The sample PCs are intrinsically linked together: how well a sample PC approximates its population counterpart is determined by the estimation error of the sample covariance matrix as well as how far its corresponding eigenvalue is from the others, and its distribution generally depends on all eigenvalues and eigenvectors of the covariance matrix. This means that inferences about one PC necessarily require information about all others and therefore impractical. Both hindrances disappear when considering ICA or PCs with multilinear structure. Figure \ref{fig:ica2} illustrates how eigengap may affect PCA but not ICA. The top left panel plots 500 observations from $U([-0.5,0.5]\times [-0.5/1.05, 0.5/1.05])$. Because the two coordinates are independent with different variances, they are both the PCs and ICs. The top middle and right panels plot the two sample PCs and sample ICs, respectively. It is clear that the sample PCs perform poorly even with such a large sample size. The bottom panels repeat the same exercise with a larger eigengap: observations are now sampled from $U([-0.5,0.5]\times [-0.5/1.5, 0.5/1.5])$. Comparing the two scenarios, the impact of eigengap on PCs is evident, whereas the ICs are virtually unaffected.
	
	\subsection{Challenges of Using Tensor Methods}
	
	While there are plenty of reasons to use tensor methods in data analysis, how to do so appropriately and take full advantage of the benefits they can offer pose many unique and fundamental challenges.
	
	Conceptually, many standard notions and intuitions we have developed for matrices need to be re-examined when extending to higher order tensors. One of the most notable examples is the best low rank approximation. In the case of matrices, the best low rank approximation can be characterized by singular value decomposition thanks to Eckart-Young theorem \citep[see e.g.,][]{bhatia2013matrix,stewart1990matrix}. In the case of higher order tensors, this is much more subtle since the best low rank approximation may not even exist! See, e.g., \cite{hackbusch2012tensor}.
	Even worse, this is not a mere pathological exception but rather a prevalent phenomenon. Consider, for example, a $2\times 2\times 2$ tensor $\scrT$ whose $(1,1,1)$ and $(2,2,2)$ entries are 1 and other entries are 0. Suppose that we observe $\scrX=\scrT+\sigma\scrE$ where $\scrE$ is an iid standard normal ensemble and want to estimate $\scrT$ knowing apriori that it has rank two. The maximum likelihood estimate (MLE) can then be \emph{defined} as the best rank-two approximation to $\scrX$. Unfortunately, it can be derived that, for any $\sigma>0$, there is a nonzero probability that the MLE defined as such does not even exist.
	
	Another ubiquitous challenge is computation. As noted by \cite{hillar2013}, ``most tensor problems are NP hard''. This means that we need to be always mindful of the computational feasibility of tensor methods when applying them to high dimensional data. Indeed one of the most intriguing phenomena that is often observed for tensor methods is that statistically optimal procedures are often found to be computationally infeasible, yet computationally tractable procedures are known to be statistically suboptimal.
	
	\subsection{Organization} 
	In this article, we shall review representative progresses that have been made to take advantage of the benefits offered by tensors while tackling the challenges they bring about in eight statistical settings. These examples are selected to illustrate how the opportunities and challenges of tensor methods manifest in different ways and under different contexts. Given the fast growing literature, the review is by no means exhaustive and the choices of examples necessarily reflect our personal and oftentimes unintentional bias.
	
	The rest of the article is organized as follows: we shall first briefly introduce some of the basic notations as well as algebraic and algorithmic tools for tensors in the next section. Section 3 discusses how relevant tensor methods have been used, what the new challenges are, how they can be addressed, and what the unique phenomena are in the eight statistical problems.
	
	\section{Background and Tools}
	
	\subsection{Notation}
	We denote vectors and matrices by boldfaced lower and upper case letters respectively: e.g., $\bv\in\RR^{d}$ and $\bM\in\RR^{d_1\times d_2}$; and tensors by script-style letters, e.g., $\scrT$. Entries of a matrix or tensor is represented by upper-case letters with a set of indices, e.g., $T_{i_1i_2\cdots i_p}$ is the $(i_1,\ldots, i_p)$ entry of an $p$th order tensor $\scrT$. We first list a set of basic notations and tensor operations that we shall use throughout the paper.
	\subsubsection{Matricization} A $p$-th order tensor $\scrT\in \RR^{d_1\times\cdots\times d_p}$ has the mode-$k$ matricization 
	$
	\calM_{k}(\scrT)\in \RR^{d_k\times \prod_{j\neq k, j\in [p]}d_j}
	$
	with elements 
	$
	\left(\calM_k(\scrT)\right)_{i_kj}
	=
	T_{i_1\dots i_{k-1}i_ki_{k+1}\dots i_p}
	$
	where 
	$$
	j=1+\sum_{l\in [p],l\neq k}(i_l-1)J_l
	\quad
	\text{with}
	\quad 
	J_l=\prod_{m\in [l-1],m\neq k}d_m.
	$$
	
	\subsubsection{Outer Product} For two tensors $\scrT\in \RR^{d_1\times\dots\times d_p}$, $\scrT'\in \RR^{d_1'\times\dots\times d_q'}$, their outer product $\scrT\circ\scrT'$ is a $(p+q)$th order tensor whose entries are given by
	$$
	(\scrT\circ\scrT')_{i_1\dots i_pj_1\dots j_q}=T_{i_1\dots i_p}T'_{j_1\dots j_q}.
	$$
	In particular, the outer product of two vectors is a matrix.
	
	\subsubsection{Multiplication} A tensor $\scrT\in \RR^{d_1\times\dots\times d_p}$ can be multiplied along its $k$-th mode by a matrix $\bM\in \RR^{m_k\times d_k}$ to get a new tensor 
	$
	\scrT'
	:=\scrT\times_k\bM\in \RR^{d_1\times d_2\times \dots d_{k-1}\times m_k\times d_{k+1}\dots \times d_p}
	$
	with elements 
	$$
	T'_{i_1\dots i_{k-1}j_ki_{k+1}\dots i_p}
	=\sum_{i_k\in [d_k]}T_{i_1\dots i_p}M_{j_ki_k}
	$$
	for $i_l\in [d_l]$ for $l=1,\dots,p$, $l\neq k$, and $j_k\in [m_k]$.

	\subsubsection{Norms} The inner product between two tensors $\scrT, \scrT'\in \RR^{d_1\times\dots\times d_p}$ is given by 
	$$
	\langle \scrT,\scrT'\rangle =\sum_{i_1,\dots,i_p}T_{i_1\dots i_p}T'_{i_1\dots i_p}.
	$$
	The Frobenius norm of a tensor is then defined as
	$
	\|\scrT\|_{\rm F}:=\langle \scrT,\scrT\rangle^{1/2}$. The spectral norm or operator norm of $\scrT$ is defined as
	$
	\|\scrT\|:=\sup_{\bu_k\in \SS^{d_k-1}}\langle \scrT, \bu_1\circ\dots\circ \bu_p\rangle,
	%
	$
	where $\SS^{d-1}$ denotes the unit sphere in $\RR^d$.
	
	\subsection{Tensor Decompositions}
	Singular value decomposition (SVD) plays a central role in analyzing data formatted in matrices. It turns out that there are multiple ways to generalize SVD to higher order tensors, most notably the so-called CP decomposition and Tucker decomposition.
	
	\subsubsection{CP Decomposition}
	One of most popular approaches to tensor decomposition is the so-called Canonical Polyadic (CP) decomposition \citep{hitchcock1927expression, carroll1970analysis, harshman1970foundations} that expresses a tensor as the sum of rank-one tensors, e.g., tensors that can be expressed as outer product of vectors. To fix ideas, consider a third order tensor $\scrT\in\RR^{d_1\times d_2\times d_3}$. Its CP decomposition has the form
	\begin{equation}\label{eq:cp-decomp}
		\scrT=\sum_{i=1}^r\lambda_i\ba_i\circ\bb_i\circ\bc_i.
	\end{equation}
	where $\lambda_i$s are scalars, and $\ba_i\in\SS^{d_1-1}$, $\bb_i\in \SS^{d_2-1}$ and $\bc_i\in\SS^{d_3-1}$ are unit vectors.

	\emph{CP Rank.} The smallest integer $r$ for which such a decomposition holds is called the CP-rank of $\scrT$. It is convenient to use the notation $\scrT=[{\boldsymbol{\lambda}};\bA,\bB,\bC]$ for CP decomposition as described by Equation~\ref{eq:cp-decomp} where ${\boldsymbol{\lambda}}=(\lambda_1,\ldots,\lambda_r)^\top$, $\bA=[\ba_1\,\ba_2\,\dots\,\ba_r]$, and similarly $\bB=[\bb_1\,\bb_2\,\dots\,\bb_r]$, and $\bC=[\bc_1\,\bc_2\,\dots\,\bc_r]$.
	
	\emph{Orthogonally Decomposable Tensors.} CP decomposition can be viewed as a generalization of SVD. In the case of matrices, the factor matrices can also be made orthonormal without loss of generality. For tensors this is generally not possible. But if a tensor does allow for a CP decomposition with orthonormal factor matrices, e.g.,
	\begin{equation*}
		\langle\ba_i,\ba_j\rangle
		=\langle\bb_i,\bb_j\rangle
		=\langle\bc_i,\bc_j\rangle
		=0,\qquad 
		\forall~ i\neq j, 
	\end{equation*}
	we say it is orthogonally decomposable (ODECO).
	
	\emph{Uniqueness.} Note that factors $\{\ba_i\}$, $\{\bb_i\}$ and $\{\bc_i\}$ in CP decomposition are not uniquely defined. For example, we can always replace $\ba_i$ by $-\ba_i$ and $\bb_i$ by $-\bb_i$, or by permuting the indices. Interestingly, however, modulo this trivial ambiguity, they may be uniquely identified. This remarkable observation was first made by J. B. Kruskal who showed \citep{kruskal1977three} that if
	\begin{equation*}
		k(\bA)+k(\bB)+k(\bC)\ge 2r-1,
	\end{equation*}
	then the CP decomposition of $\scrT$ given by Equation~\ref{eq:cp-decomp} is unique up to a permutation of the indices $\{1,\dots,r\}$ and sign changes of the factors. Here for a matrix $\bM\in \RR^{d\times r}$, $k(\bM)$ denotes its Kruskal rank, that is the maximum integer $k$ such that any $k$ columns of $\bM$ are linearly independent. This immediately implies that ODECO tensors are identifiable.
	
	\subsubsection{Tucker Decomposition}
	Tucker decomposition \citep{tucker1963} is another popular way to generalize SVD to tensors. More specifically, the Tucker decomposition of a third order tensor $\scrT\in\RR^{d_1\times d_2\times d_3}$ can be written as:
	\begin{equation*}
		\scrT=[\scrC;\bU_1,\bU_2,\bU_3]:=\scrC\times_1\bU_1\times_2\bU_2\times_3\bU_3,
	\end{equation*}
	where $\scrC\in \RR^{r_1\times r_2\times r_3}$ is a so-called core tensor, and $\bU_q\in \RR^{d_q\times r_q}$ is the $q$-th mode component matrices, for $q=1,2,3$. Note that Tucker decomposition is not unique. For example, for any invertible matrix $\bM\in \RR^{r_1\times r_1}$, 
	$$\scrT=[\scrC;\bU_1,\bU_2,\bU_3]=[\scrC\times_1\bM;\bU_1\bM^{-\top},\bU_2,\bU_3].$$ 
	
	\emph{Higher Order SVD.} Without loss of generality, we can always take component matrices $\bU^{(q)}$s to be orthonormal as a generalization of the matrix SVD. Note that
	\begin{equation*}
		\calM_q(\scrT)=\bU_q\calM_q(\scrC)(\otimes_{q'\neq q}\bU_{q'})^{\top}
	\end{equation*}
	for $q=1,2,3$ where $\otimes$ stands for Kronecker product. \cite{de2000multilinear} showed that we can force the rows of $\calM_q(\scrC)$ to be orthogonal for all $q=1,2,3$ so that $\bU_q\in \RR^{d_q\times r_q}$ is the left singular matrix of $\calM_q(\scrT)$, and the $\ell_2$ norm of the rows of $\scrM_q(\scrC)$ are the singular values. This is the so-called higher order SVD (HOSVD).
	
	\emph{Multilinear ranks.} Tucker decomposition also corresponds to a new notion of rank for tensors. In particular, we can take $r_q$ to be the rank of $\calM_q(\scrT)$ and the $p$-tuple $(r_1,\dots,r_p)$ is often referred to as the multilinear rank or Tucker rank of $\scrT$. It is worth noting that, in general, $r_q$s are different.

	\subsection{Algorithms}
	
	In a nutshell, most of the tensor methods for data analysis can be reduced to finding a low-rank, either CP rank or Tucker rank, tensor that fits the data well. This is often cast as minimizing a loss function $L(\cdot)$ over tensors satisfying certain rank constraints. Such problems are generally highly nonconvex and computationally challenging. We now discuss several popular practical strategies.
	
	\subsubsection{Alternating Minimization} The idea behind alternating minimization is that we can reparametrize a low rank tensor by its decomposition and then optimize over the factors in an iterative fashion. Consider, for example, optimizing $L$ over $\calF(d_1,d_2,d_3,r_1,r_2,r_3)$,  which is a certain subset of tensors in $\RR^{d_1\times d_2\times d_3}$ with multilinear ranks $(r_1,r_2,r_3)$.  Additional structural information such as sparsity or incoherence may also be incorporated in $\calF(d_1,d_2,d_3,r_1,r_2,r_3)$.   The alternating minimization algorithm operates on the joint space induced by the tensor factorization form. 
	By re-parametrizing a tensor $\scrT \in \calF(d_1,d_2,d_3,r_1,r_2,r_3)$ using Tucker decomposition,  the optimization problem can be equivalently written as 
	\begin{equation*}
		\min \ L\big([\scrC; \bU_1, \bU_2, \bU_3]\big), \quad {\rm s.t.}\quad \scrC\in\RR^{r_1\times r_2\times r_3},\ \bU_k\in \calU(d_k, r_k),\ k=1,2,3.
	\end{equation*}
	The feasible set $\calU(d_k, r_k)$ is a subset of $\RR^{d_k\times r_k}$, predetermined according to particular applications or for identifiability purposes. The alternating minimization algorithm is an iterative procedure that minimizes the loss function by alternately optimizing over individual parameters. The pseudocodes for alternating minimization are provided in Algorithm~\ref{alg:alt-min}. The alternating minimization algorithm usually converges fast provided that a sufficiently good initialization is given and the sub-problems in Equation~\ref{eq:alt-min-update} can be solved efficiently.
	
	\begin{algorithm}
		\caption{Alternating Minimization}\label{alg:alt-min}
		\begin{algorithmic}
			\STATE{\textbf{Input}: loss $L\big([\scrC; \bU_1, \bU_2, \bU_3]\big)$; feasible sets $\bU_k\in\calU(d_k, r_k)\subset \RR^{d_k\times r_k}$ for $\forall k=1,2,3$ }
			\STATE{Initialization: $\hat{\scrC}^{[0]}\in\RR^{r_1\times r_2\times r_3}$ and $\hat\bU_k^{[0]}\in\calU(d_k, r_k)$ for $\forall k=1,2,3$; }
			\STATE{Set $t\leftarrow 0$.}
			\WHILE{{\it not converged}}
			\STATE{Update by }
			\begin{align}
				\hat\bU_1^{[t+1]} &\longleftarrow \underset{\bU\in\calU(d_1, r_1)}{\arg\min} \ L\Big(\big[\hat\scrC^{[t]};  \bU, \hat \bU_2^{[t]}, \hat \bU_3^{[t]}\big]\Big)\notag\\
				\hat\bU_2^{[t+1]} &\longleftarrow \underset{\bU\in\calU(d_2, r_2)}{\arg\min} \ L\Big(\big[\hat\scrC^{[t]}; \hat \bU_1^{[t+1]}, \bU, \hat \bU_3^{[t]}\big]\Big)\label{eq:alt-min-update}\\
				\hat\bU_3^{[t+1]} &\longleftarrow \underset{\bU\in\calU(d_3, r_3)}{\arg\min} \ L\Big(\big[\hat\scrC^{[t]}; \hat \bU_1^{[t+1]}, \hat \bU_2^{[t+1]}, \bU\big]\Big)\notag\\
				\hat\scrC^{[t+1]} &\longleftarrow \underset{\scrC\in\RR^{r_1\times r_2\times r_3}}{\arg\min} \ L\Big(\big[\scrC; \hat \bU_1^{[t+1]}, \hat \bU_2^{[t+1]}, \hat \bU_3^{[t+1]}\big]\Big)\notag
			\end{align}				
			Set $t\leftarrow t+1$.		
			\ENDWHILE
			\STATE{\textbf{Return}: $\hat \scrT^{[t]}=\big[\hat \scrC^{[t]}; \hat\bU_1^{[t]},\hat\bU^{[t]}_{2}, \hat\bU_3^{[t]}\big]$.}
		\end{algorithmic}
	\end{algorithm}

	\emph{Higher-order orthogonal iteration.} The higher-order orthogonal iteration (HOOI) is a special case of alternating minimization algorithm applied to the loss function $L(\scrT)=\|\scrX-\scrT\|_{\rm F}^2$ with the feasible set $\calU(d_k, r_k)=\OO_{d_k,r_k}$, the collection of $d_k\times r_k$ matrices with orthonormal columns. Here, $\scrX$ represents a $d_1\times d_2\times d_3$ data tensor. The sub-problems Equation~\ref{eq:alt-min-update} admit an explicit solution, for example, 
	$$
	\hat \bU^{[t+1]}_{1}={\rm SVD}_{r_1}\Big(\calM_1\big(\scrX\times_2 \hat\bU^{[t]\top}_{2}\times_3 \hat\bU^{[t]\top}_{3}\big)\Big),
	$$
	where ${\rm SVD}_r(\cdot)$ returns the top-$r$ left singular vectors of a matrix. HOOI was initially introduced by \cite{de2000best} to study the best low-rank approximation of a given tensor. Variants of HOOI have been designed to accommodate additional structures. See, e.g., \cite{allen2012sparse,sun2017provable,zhang2019optimal, xia2021statistically, ke2019community, jing2021community}. 
	
	\emph{Power Iteration.} Tensor power iteration (tens-PI) is a special case of HOOI used for finding the best rank-one approximation of a given tensor. It most often appears in alternating minimization for optimizing $L$ over the set of CP-decomposable tensors. By re-parametrizing a tensor  $\scrT$ using CP-decomposition,  we aim to minimize $L\big(\sum_{i=1}^r \lambda_i \ba_i\circ \bb_i\circ \bc_i\big)$,  where $\lambda_i$s are scalars and $\ba_i, \bb_i, \bc_i$s are unit vectors. The algorithm alternately minimizes the loss function over individual rank-one components. This approach leads to the tensor power iteration algorithm if $L(\scrT)=\|\scrX-\scrT\|_{\rm F}^2$ for a given data tensor $\scrX$. For instance, given the estimates $\big\{\hat\lambda_i^{[t]}, \hat\ba_i^{[t]}, \hat\bb_i^{[t]}, \hat\bc_i^{[t]}\big\}_{i=1}^r$ at the $t$-th iteration, the alternating minimization algorithm updates by solving 
	\begin{align*}\label{eq:power-iter-update}
		(\hat\lambda_j^{[t+1]}, \hat\ba_j^{[t+1]}, \hat\bb_j^{[t+1]}, \hat\bc_j^{[t+1]})\longleftarrow \underset{\lambda, \ba, \bb, \bc}{\arg\min}\ \Big\|\Big(\scrX-\sum\nolimits_{i=1}^{j-1} \hat\lambda_i^{[t+1]} \hat \ba_i^{[t+1]}\circ \hat \bb_i^{[t+1]}\circ \hat \bc_i^{[t+1]}\\\nonumber-\sum\nolimits_{i=j+1}^{r} \hat\lambda_i^{[t]} \hat \ba_i^{[t]}\circ \hat \bb_i^{[t]}\circ \hat \bc_i^{[t]}\Big)-\lambda \ba\circ \bb\circ \bc \Big\|_{\rm F}^2,
	\end{align*}
	which is often solved by the tensor power iteration algorithm. Tensor power iteration has been used in many works including \cite{anandkumar2014tensor, anandkumar2014guaranteed,mu2015successive,mu2017greedy,sharan2017orthogonalized}. In statistical applications, tens-PI is used with context aware smart initializations for SVD of the covariance tensor in \cite{liu2022characterizing}, for learning mixture models in \cite{anand2014sample}, and for ICA in \cite{auddy2023large}. A truncated tens-PI has been used for sparse CP decomposition in \cite{sun2017provable}. Tens-PI has been proven effective for tensor completion in \cite{jain2014provable}.

	\subsubsection{Gradient Descent} Another common strategy for optimizing over low rank tensors is gradient descent.  To minimize $L$ over $\calF(d_1,d_2,d_3, r_1, r_2, r_3)$,  gradient descent algorithms can operate  directly on the space of $d_1\times d_2\times d_3$ tensors with multilinear ranks $(r_1,r_2,r_3)$.  
	Multiple choices of gradients are available,  such as the standard (vanilla) gradient in Euclidean space and the Riemannian gradient in the associated manifolds.  The general architect of gradient descent algorithms can be described by Algorithm~\ref{alg:grad-desc}.  Here,  $\textsf{Proj}_{\calF(d_1,d_2,d_3,r_1,r_2,r_3)}(\cdot)$ denotes the projection operator, which enforces that each update resides in $\calF(d_1,d_2,d_3,r_1,r_2,r_3)$.

	\begin{algorithm}
		\caption{Gradient Descent}\label{alg:grad-desc}
		\begin{algorithmic}
			\STATE{\textbf{Input}: loss $L(\scrT)$;  initializations $\hat \scrT^{[0]}\in\calF(d_1,d_2,d_3,r_1,r_2,r_3)$; step size $\eta>0$;  $t=0$ }
			\WHILE{{\it not converged}}
			\STATE{Compute a \emph{certain} type of gradient: $\hat\scrG^{[t]}\longleftarrow \textsf{grad}\big( L(\hat\scrT^{[t]})\big)$}
			\STATE{Update by
				$
				\hat \scrT^{[t+1]}\longleftarrow \textsf{Proj}_{\calF(d_1,d_2,d_3,r_1,r_2,r_3)}\Big(\hat \scrT^{[t]}-\eta\cdot \hat\scrG^{[t]}\Big)
				$}
			\STATE{Set $t\leftarrow t+1$}			
			\ENDWHILE
			\STATE{\textbf{Return}: $\hat \scrT^{[t]}$.}
		\end{algorithmic}
	\end{algorithm}

	\emph{Projected Gradient Descent.} When the vanilla gradient is computed with $\hat\scrG^{[t]}:=\nabla L(\hat\scrT^{[t]})$,   it results in the projected gradient descent algorithm.  The gradient $\hat\scrG^{[t]}$ is typically a full-rank tensor, as is the gradient descent update  $\hat\scrT^{[t]}-\eta\cdot \hat\scrG^{[t]}$.  The projection step entails finding a low-rank approximation of the given full-rank tensor.  Although finding the best low-rank approximation is a challenge for tensors,  a reasonably good approximation is sufficient to ensure the convergence of Algorithm~\ref{alg:grad-desc}.  As demonstrated in \cite{chen2019non}, HOSVD can serve as the projection step, and Algorithm~\ref{alg:grad-desc} converges at a linear rate,  yielding a statistically optimal estimator under various models.

	\emph{Riemannian Gradient Descent.} By viewing $\calF(d_1,d_2, d_3, r_1,r_2,r_3)$ as a Riemannian manifold,  one can take $\hat\scrG^{[t]}$ as the Riemannian gradient,  which leads to the Riemannian gradient descent algorithm \citep{edelman1998geometry,  kressner2014low}.  The Riemannian gradient has multilinear ranks at most $(2r_1, 2r_2, 2r_3)$ ensuring that the gradient descent update has multilinear ranks at most $(3r_1, 3r_2, 3r_3)$.  This low-rank structure enables much faster computation for the subsequent projection step.  Riemannian gradient descent was investigated by \cite{cai2022generalized} for robust tensor decomposition and by \cite{cai2022provable} for  noisy tensor completion.

	\emph{Grassmannian and Joint Gradient Descent.} Gradient descent algorithms can also operate on the joint space induced by the tensor factorization form.  To minimize  $L\big([\scrC; \bU_1, \bU_2,  \bU_3]\big)$,  one can compute the gradients with respect to $\scrC$ and $\bU_{k}$s,  respectively.  The constraint $\bU_{k}\in\OO_{d_k,r_k}$ is often imposed for identifiability purposes.  By viewing $\OO_{d_k,r_k}$ as the Grassmannian manifold,  we can utilize the Grassmannian gradient to update $\hat\bU_{k}$,  leading to the Grassmannian gradient descent algorithm \citep{xia2019polynomial,  lyu2023latent}.   Instead of enforcing orthogonality constraints,  \cite{han2022optimal} considers minimizing $L\big([\scrC; \bU_1, \bU_2,  \bU_3]\big)$ with the penalty $\lambda_1 \sum_{k=1}^3 \|\bU^{\top}_k\bU_{k}-\bI_{r_3}\|_{\rm F}^2$ and proposes a vanilla gradient descent algorithm,  jointly with respect to $\scrC$ and $\bU_{k}$s.  
	
	\subsubsection{Overcoming Nonconvexity} Problems involving low rank tensors are usually highly nonconvex with many local optima. Either alternating minimization or gradient descent algorithms are prone to be stuck with bad local optima. An intriguing and recurring phenomenon that is often observed is that with a good initialization, one can usually obtain a statistically efficient estimate with an optimal signal-to-noise ratio. Yet to overcome the nonconvexity, a much higher signal-to-noise ratio is necessary. Nonetheless, several strategies are often used to alleviate the challenge.
	
	
	\emph{Spectral Initialization.} Depending on the data generating mechanism, different types of spectral initialization are commonly used. In particular, HOSVD was widely used for initialization when considering tensors of low Tucker ranks. See, e.g., \cite{zhang2018tensor, cai2019nonconvex, montanari2018spectral, ke2019community}. Both matricization, e.g., unfolding a tensor into matrices, and random slicing, e.g., random linear combinations of all slices, are often used for ODECO tensors and more generally in the context of CP decomposition. See, e.g., \cite{anandkumar2014tensor,anand2014sample,auddy2023large}. \cite{xia2019polynomial} introduced a spectral initialization approach based on a second-order U-statistics that has been shown to be effective in a number of settings including tensor completion and regression \citep[see, e.g.,][]{cai2022provable, xia2021statistically}.

	\emph{Convex Relaxations.} Another popular strategy to overcome the nonconvexity is by convex relaxation. This is usually accomplished by matricization and resort to matrix-based convex approaches. See, e.g., \cite{gandy2011tensor, mu2014square,raskutti2019convex}. More powerful semidefinite or sum-of-squares relaxation approaches are also commonly considered. See, e.g., \cite{jiang2015tensor,tomioka2013convex,hopkins2015tensor,ma2016polynomial,schramm2017fast}. 
	
	\emph{Other Approaches.} One common way to avoid being stuck at a local optimum is by adding occasional and carefully designed jumps to the original iterative algorithms. See, e.g., \cite{belkin2018eigenvectors}. Adding stochastic noise to gradients is another way to escape the trap of local optima  \citep[see, e.g.,][]{ge2015escaping}. The approximate message passing algorithm has been proven successful in learning rank-one tensor PCA, as demonstrated by \cite{richard2014statistical}. 
	
	\section{Statistical Models}
	We now consider several specific tensor-related statistical problems. To fix ideas, unless otherwise indicated, we shall focus our discussion on third order cubic tensors, e.g., $\scrT\in \RR^{d\times d\times d}$. Oftentimes, more general results are available but we opt for this simplification for brevity.
	
	\subsection{Tensor SVD}
	Let us begin with a canonical example for studying tensor methods where we wish to recover the decomposition of a signal tensor $\scrT$ from an observation contaminated with noise, e.g., $\scrX=\scrT+\scrE$. For simplicity, we assume that the noise $\scrE$ is an iid ensemble, e.g., its entries are independently sampled from a common distribution.
	
	\subsubsection{Computational Gap} Tensor SVD provides a simple statistical model to understand the role of computational considerations in tenor related problems.
	
	\emph{Gaussian Noise.} The simplest case is when the signal tensor $\scrT$ has rank one, e.g., $\scrT=[\lambda; \ba, \bb,\bc]=\lambda\ba \circ\bb\circ\bc$ where $\lambda>0$, $\ba,\bb,\bc\in \SS^{d-1}$. When $\scrE$ has independent standard Gaussian entries, the maximum likelihood estimate (MLE) is given by 
	\begin{equation}\label{eq:rank1-tens-svd}
		(\hat{\lambda},\hat{\ba}_{\rm SVD},\hat{\bb}_{\rm SVD},\hat{\bc}_{\rm SVD})
		:=\underset{\gamma>0,\bx,\by,\bz\in\SS^{d-1}}
		{\arg\min}
		\|\scrX-\gamma\bx\circ\by\circ\bz\|_{\rm F}^2.
	\end{equation}
	\cite{richard2014statistical} showed that for $\lambda\gtrsim\sqrt{d}$, the MLE is indeed consistent, and conversely, if $\lambda\lesssim\sqrt{d}$, then no consistent estimates of $\ba, \bb, \bc$ exist. 
	
	However, as proved by \cite{arous2019landscape}, the optimization problem on the right hand side of Equation \ref{eq:rank1-tens-svd} can have exponentially, in $d$, many spurious local minima when $\lambda\ll d$. This computational challenge can be overcome by, for example, spectral initialization after matricization. More specifically, the right hand side of Equation \ref{eq:rank1-tens-svd} can be optimized by power iteration intialized by HOSVD, e.g., $\ba^{[0]}$, $\bb^{[0]}$ and $\bc^{[0]}$ are the leading left singular vector of $\calM_1(\scrX)$, $\calM_2(\scrX)$ and $\calM_3(\scrX)$ respectively. This approach yields consistent estimate of $\scrT$ when $\lambda\gtrsim d^{3/4}$. The same signal-to-noise ratio requirement can also be achieved by a number of alternative algorithms including approximate message passing \citep[see, e.g.,][]{richard2014statistical}. Nonetheless, there remains an additional cost in the signal-to-noise ratio required for consistent estimate for computationally tractable methods.
	
\begin{figure}[htbp]
		\centering
		\includegraphics{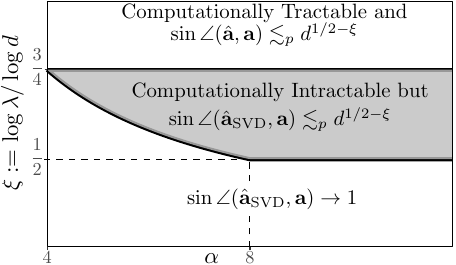}
		\caption{Statistical and computational tradeoff for tensor SVD: horizontal axis represents the moment condition for noise, vertical axis corresponds to the signal strength. Reproduced and modified from \cite{auddy2021estimating}.}\label{fig:rank1-svd}
	\end{figure}

	\emph{General Noise.} This observation goes beyond consistency and Gaussian noise. As shown by \cite{auddy2021estimating}, the gap between statistical efficiency and computational tractability can be characterized by moment conditions of the entries of $\scrE$. Specifically, consider the case $\lambda=d^\xi$ and $\EE|E_{111}|^\alpha<\infty$ for some $\alpha\ge 4$. Then the minimax optimal rate of convergence is given by
	$$
	\inf_{\tilde{\ba},\tilde{\bb},\tilde{\bc}}\sup_{\ba,\bb,\bc\in \SS^{d-1}}\max\{\sin\angle (\tilde{\ba},\ba), \sin\angle (\tilde{\bb},\bb),\sin\angle (\tilde{\bc},\bc)\}\asymp d^{1/2-\xi},
	$$
	where the infimum is taken over all estimates based on $\scrX$. This rate of convergence can be attained, in particular, by $\hat{\ba}_{\rm SVD},\hat{\bb}_{\rm SVD},\hat{\bc}_{\rm SVD}$ when $\xi>\max\{1/2,1/4+2/\alpha\}$. On the other hand, computationally tractable algorithms they developed can only reach minimax optimality when $\xi>3/4$, thus creating a gap in signal-to-noise ratio requirement. See Figure \ref{fig:rank1-svd}. It is interesting to note that for any $\alpha>8$, the gap is identical to the Gaussian case but narrows as $\alpha$ decreases to 4.
	
	\emph{Low Tucker Rank Spikes.} Generalizing beyond rank one case in a Tucker decomposition framework, one may consider $\scrT=[\scrC; \bU_1,\bU_2,\bU_3]$
	%
	where $\scrC\in\RR^{r\times r\times r}$ is the core tensor and we wish to estimate the semi-orthonormal component matrices $\bU_q\in \RR^{d\times r}$ for $q=1,2,3$. Under the Gaussian noise, \cite{zhang2018tensor} showed that similar computational gap exists. More specifically, if $r\lesssim\sqrt{d}$ and $\lambda_{\min}:=\min_{q=1,2,3}\sigma_{r}(\calM_q(\scrC))\gtrsim\lambda_{\rm comp}\asymp d^{3/4}$, then
	$$
	\sin\Theta\left(\hat{\bU}_q,\bU_q\right)\lesssim \sqrt{d}/\lambda_{\min}\quad 
	\text{ for }q=1,2,3,
	$$
	where $\hat{\bU}_q$ is the estimate computed by HOOI initialized via HOSVD. Furthermore, this rate is also minimax optimal.

	\subsubsection{Perturbation Bounds}\label{sec:perturb}
	As in the matrix case, perturbation bounds play a crucial role in the analysis of low rank tensors in general and tensor SVD in particular. For ODECO tensors, deterministic perturbation bounds for singular values and vectors, in a spirit similar to classical results for matrices such as those due to Weyl, Davis, Kahan and Wedin, can be established. More specifically, assume that $\scrT=\sum_{i=1}^d \lambda_i\ba_i\circ\bb_i\circ\bc_i$ is an ODECO tensor. Let $\hat{\scrT}=\sum_{i=1}^d \hat{\lambda}_i\hat{\ba}_i\circ\hat{\bb}_i\circ\hat{\bc}_i$ be an (arbitrary) estimate of $\scrT$ that is also ODECO. \cite{auddy2020perturbation} showed that there exists a numerical constant $C>1$ such that
	$$
	\max_{1\le i\le d}|\hat{\lambda}_i-\lambda_i|\le C \|\hat{\scrT}-\scrT\|
	$$
	and there is a permutation $\pi:[d]\to[d]$ such that
	$$
	\max \{\sin\angle(\hat{\ba}_{\pi(i)},\ba_i),\sin\angle(\hat{\bb}_{\pi(i)},\bb_i),\sin\angle(\hat{\bc}_{\pi(i)},\bc_i)\}\le {C\|\hat{\scrT}-\scrT\|\over \lambda_i}.
	$$
	Although these bounds are similar in spirit to those for matrices, there are also crucial distinctions. In particular, the $\sin\Theta$ theorems of Davis-Kahan-Wedin bound the perturbation effect on the $i$th singular vector by $C\|\hat{\scrT}-\scrT\|/\min_{i'\neq i} |\lambda_{i'}-\lambda_i|$. The dependence on the gap between $\lambda_i$ and other singular values is unavoidable for matrices. This is not the case for higher-order ODECO tensors where perturbation affects the singular vectors in separation. Note that the permutation $\pi$ of indices is necessary because we do not require that $\lambda_i$s are distinct or sufficiently apart from each other. 
	
	Furthermore, if the estimation error is small when compared with signal, e.g., $\|\hat{\scrT}-\scrT\|=o(\lambda_i)$, then we can take the constant in the above bounds to be one:
	\begin{equation*}
		|\hat{\lambda}_{\pi(i)}-\lambda_i|\le \|\hat{\scrT}-\scrT\|    
	\end{equation*}
	and
	\begin{equation}\label{eq:odec-pert2}
		\max \left\{\sin\angle(\hat{\ba}_{\pi(i)},\ba_i),\sin\angle(\hat{\bb}_{\pi(i)},\bb_i),\sin\angle(\hat{\bc}_{\pi(i)},\bc_i)\right\}
		\le {\|\hat{\scrT}-\scrT\|\over \lambda_i}(1+o(1)).  
	\end{equation}
	Both bounds can be shown to be sharply optimal. See \cite{auddy2020perturbation} for further detailed discussion. 
	
	Applying these perturbation bounds one can see that the minimax optimal rates for estimating ODECO $\scrT$ are essentially the same as those for rank-one spikes and can be attained, in particular, by $\hat{\scrT}=\argmin_{\scrA\ {\rm is\ ODECO}}\|\scrX-\scrA\|$. This estimate, as in the rank-one spike case, is not amenable for computation. There is a fruitful line of recent research in developing computationally tractable estimates and algorithm dependent bounds, which again incates a gap in performance between computationlly tractable and intractable methods. See, e.g., \cite{anandkumar2014tensor,luo2021sharp,janzamin2019spectral}.

	
	\subsection{Multiway PCA}
	Tensor SVD is closely related to multiway PCA. Traditional PCA is a useful tool for dimension reduction but it does not account for the multiway structure when each observation consists of a matrix or a possibly higher order array. Naively applying PCA to multiway observations by ``flattening'' them into vectors is statistically inefficient and also lacks interpretability. To fix ideas, let us focus here on the case when each observation is an array of dimension $d\times d\times d$.
	
	Recall that, for a random vector $\bX\in \RR^d$, the first PC direction is the vector $\ba\in \SS^{d-1}$ which maximizes  $\sigma^2(\ba):={\rm Var}(\langle\bX,\ba\rangle)$. Similarly, for a tensor $\scrX\in\RR^{d\times d\times d}$, it is natural to define the leading multiway principal component $\scrU_1\in\RR^{d\times d\times d}$ as:
	\begin{equation*}
		\scrU_1=
		\underset{\ba\circ\bb\circ \bc,\ \ba,\bb,\bc\in\SS^{d-1}}
		{\arg\max}{\rm Var}(\langle\scrX,\ba\circ\bb\circ \bc\rangle).
	\end{equation*}
	The constraint that $\scrU_1$ is rank one, i.e., $\scrU_1=\bu_1^{(1)}\circ\bu_1^{(2)}\circ \bu_1^{(3)}$ for $\bu_1^{(1)}, \bu_1^{(2)},\bu_1^{(3)}\in\SS^{d-1}$, helps ensure that $\scrU_1$ conforms to the multiway nature of $\scrX$. Other PCs can be defined successively:
	\begin{equation*}
		\scrU_k=\underset{\substack{
				\ba\circ\bb\circ \bc,\ 
				\ba, \bb, \bc\in\SS^{d-1},\\
				\ba\perp\ba_{l},\bb\perp \bb_{l},\bc\perp \bc_{l}, l<k}}{\arg\max}
		{\rm Var}(\langle\scrX,\ba\circ\bb\circ \bc\rangle).
	\end{equation*}
	Once again, similar to classical PCA, the successive PCs are defined to maximize variance in the directions orthogonal to previously found PCs $\scrU_1,\dots,\scrU_{k-1}$. The sample multiway PCs can be defined similarly with population variance now replaced by sample variances.
	
	\subsubsection{PCA or SVD} When each observation is a vector, sample PCs can be identified with SVD of the data matrix where each row corresponds to an observation, after appropriate centering. See, e.g., \cite{jolliffe2002principal}. As such, one often uses the two terms, PCA and SVD, interchangeably. The equivalence, interestingly, does not hold for multiway PCA in general even though the subtle differences between the two are often mistakenly neglected.
	
	It is easy to see that
	\begin{equation}\label{eq:multpca}\sigma^2(\scrU_1)\scrU_1\circ \scrU_1+\cdots+\sigma^2(\scrU_r)\scrU_r\circ \scrU_r
	\end{equation}
	is a rank-$r$ greedy ODECO approximation to the covariance operator $\mathscr{T}_{\Sigma}={\rm Cov}(\scrX)$,  which is a $d\times d\dots\times d$ sixth order array with entries $\sigma_{i_1i_2i_3j_1j_2j_3}={\rm Cov}(X_{i_1i_2i_3},X_{j_1j_2j_3})$. Likewise, the sample counterpart of Equation \ref{eq:multpca} is a rank-$r$ greedy ODECO approximation to the sample covariance operator. In general, however, $\scrU_k$s cannot be identified with best low rank approximations to the data tensor $\scrX_n\in \RR^{d\times d\times d\times n}$ where each mode-four slice is one observation.

	\subsubsection{Asymptotic Independence} \cite{ouyang2023multiway} showed that there are many statistical benefits to consider multiway PCA. Consider a spiked covariance model: 
	$$
	\scrX=\sum_{k=1}^d\sigma_k\theta_k\scrU_k+\sigma_0\scrE,
	$$
	where $(\theta_1,\dots,\theta_r)^{\top}\sim N(0,\II_r) $ are the random factors, and $\scrU_k=\bu_k^{(1)}\circ\bu_k^{(2)}\circ \bu_k^{(3)}\in\RR^{d\times d\times d}$ are rank one principal components with $\langle \bu_k^{(q)},\bu_l^{(q)}\rangle=1$ if $k=l$ and $0$ otherwise, for $q=1,2,3$. Moreover $\scrE$ is a noise tensor with independent $N(0,1)$ entries. Under this model,
	$$
	\sin\angle\left(\hat{\scrU}_k,\scrU_{\pi(k)}\right)\le
	C 
	\left(
	\dfrac{\sigma_0}{\sigma_{\pi(k)}}+\dfrac{\sigma_0^2}{\sigma_{\pi(k)}^2}
	\right)
	\max\left\{\dfrac{d}{n},\sqrt{\dfrac{d}{n}}\right\}
	$$
	for some permutation $\pi:[d]\to[d]$. This highlights two main benefits of using multiway PCs instead of using classical PCA after stringing the observations as long vectors. Firstly, the dimension dependence is only through $d$, which is smaller than the real dimension of $\scrU_k$, i.e., $d^3$. Secondly, the result does not require any eigengap, and allows for repeated $\sigma_k$ values. The issue of repeated singular values however leads to bias, since in this case it is difficult to identify which sample PC estimates which $\scrU_k$. 
	This problem can be overcome by a sample splitting scheme as described by \cite{ouyang2023multiway}, leading to an asymptotically normal estimate $\hat{\scrU}_k$ of $\scrU_k$. In particular, in the fixed dimension case,
	$$
	\sqrt{n}
	\left({\rm vec}(\hat{\bU}^{(q)})-{\rm vec}(\bU^{(q)})\right)
	\stackrel{d}{\to}
	N\left(
	0,{\rm diag}(\Gamma_1,\dots,\Gamma_d)
	\right)
	\text{ as }
	n\to\infty,
	$$
	where $\Gamma_k=\left(\sigma_0^2/\sigma_k^2+\sigma_0^4/\sigma_k^4\right)
	\left(\II_{d}-\bu_k^{(q)}(\bu_k^{(q)})^{\top}\right)$. Here $\hat{\bU}^{(q)}=[\tilde{\bu}_1^{(q)}\,\tilde{\bu}_2^{(q)}\,\dots\,\tilde{\bu}_d^{(q)}]$ and $\bU^{(q)}=[\bu_1^{(q)}\,\bu_2^{(q)}\,\dots\,\bu_d^{(q)}]$. Notice that the asymptotic distribution of the bias corrected sample multiway PCs are normal distributions centered at the true multiway PCs. Perhaps more importantly, the sample multiway PCs are asymptotically independent. Similar properties continue to hold in the case of growing dimensions. 
	
	\subsubsection{Principal Subspaces} One can also generalize the notion of principal subspaces to multiway data. For example, \cite{tang2023mode} introduced a spiked covariance model where
	$$
	{\rm Cov}(\scrX) = \scrT_{\Sigma_0} + \sigma^2 {\scrI_{d^3}},
	$$
	where 
	$$
	\scrT_{\Sigma_0} \times_{k = 1}^3 \bU_{k\perp} = 0,
	$$
	for some $\bU_k\in \mathbb{O}_{d, r}$, and $\scrI_{d^3}$ is the sixth order tensor with entries $I_{j_1j_2\dots j_6}=1$ if $j_1=j_4,j_2=j_5,j_3=j_6$ and $I_{j_1j_2\dots j_6}=0$ otherwise. In other words, the principal subspaces in each mode can be characterizied by orthonormal matrices $\bU_1,\bU_2,\bU_3\in \RR^{d\times r}$. In particular, they show that this is equivalent to
	$$
	\scrX = \scrM+\scrA_1 \times_1 \bU_1 +\scrA_2 \times_2 \bU_2+ \scrA_3 \times_3 \bU_3+ \scrZ,
	$$
	where $\scrM$ is a fixed mean, $\scrA_1 \in \mathbb{R}^{r\times d\times d}$, $\scrA_2 \in \mathbb{R}^{d\times r\times d}$ and $\scrA_3 \in \mathbb{R}^{d\times d\times r}$ are random tensors with mean 0, and $\scrZ \in \mathbb{R}^{d\times d\times d}$ is a noise tensor with mean 0, and covariance $\sigma^2 \scrI_{d^3}$, and is uncorrelated with random tensors $\scrA_1,\scrA_2, \scrA_3$.
	
	\subsection{Independent Component Analysis} 
	ICA is a popular method for separating data streams into independent source components and has been extensively studied. In the ICA model, one observes a random vector $\bX\in\RR^d$ where $\bX$ can be written as a linear combination of independent sources. That is, there is a mixing matrix $\bA\in \RR^{d\times d}$ such that $\bX=\bA\bS$, where $\bS=(S_1,\dots,S_d)^{\top}$ is a vector of independent random variables $S_k$ for $k=1,\dots,d$. Given $n$ independent copies of $\bX$, the goal in ICA is to recover the unknown matrix $\bA$. The problem is intimately related to tensor decomposition.
	
	To fix ideas, suppose that $\bX$ are suitably centered and scaled, i.e., $\EE(\bX)=\EE(\bS)=\mathbf{0}$ and ${\rm Var}(\bX)={\rm Var}(\bS)=\II_d$. In practice, this is often referred to as \textit{pre-whitening}. In doing so, $\bA$ becomes an orthonormal rotation matrix. The rotation matrix $\bA$ can be identifiable from observing $\bX=\bA\bS$ if and only if all but one source $S_k$ are non-Gaussian \citep{comon1992independent}. ICA methods have been developed to exploit the \textit{non-Gaussianity} of $\bX$, particularly through a nonzero kurtosis. The kurtosis of $\bX$ sampled from a prewhitened ICA model has the form \cite[see, e.g.,][]{comon2010handbook}:
	\begin{equation}\label{eq:ica-kurt}
		\mathscr{K}_4(\bX)
		=\EE(\bX\circ\bX\circ\bX\circ \bX) -\scrM_0
		=\sum_{k=1}^d
		\kappa_4(S_k)\ba_k\circ\ba_k\circ\ba_k\circ\ba_k,
	\end{equation}
	where  $
	\scrM_0=
	\EE(\bZ\circ\bZ\circ\bZ\circ\bZ)
	$
	for $\bZ\sim N(\mathbf{0},\II_d)$, 
	is a known 4-th order tensor. Moreover, $\kappa_4(S_k)=\EE S_k^4-3$ is the excess kurtosis of the $k$-th independent source $S_k$.

	\subsubsection{Perturbation Bounds} In light of Equation~\ref{eq:ica-kurt}, we can derive an estimate of $\bA$ from that of $\mathscr{K}_4(\bX)$. Perhaps the most natural choice is the sample cumulant tensor:
	$$
	\hat{\scrK}^{\rm sample}_4:=
	\dfrac{1}{n}
	\sum_{i=1}^n\bX_i\circ\bX_i\circ\bX_i\circ\bX_i
	-\scrM_0.
	$$
	As noted in Section \ref{sec:perturb}, we can then proceed to estimate $\ba_k$s by the eigenvectors, denoted by $\hat{\ba}_k$s, of an ODECO approximation of $\hat{\scrK}^{\rm sample}_4$. Using the perturbation bounds from \cite{auddy2020perturbation}, we immediately get
	$$
	\sin\angle(\hat{\ba}_{\pi(k)},\ba_k)\lesssim {\|\hat{\scrK}^{\rm sample}_4-\mathscr{K}_4(\bX)\|\over \kappa_4(S_k)}.
	$$
	Interestingly though, $\hat{\scrK}^{\rm sample}_4$ is inconsistent in that $\|\hat{\scrK}^{\rm sample}_4-\mathscr{K}_4(\bX)\|\not\to_p 0$ unless $n\gtrsim d^2$. Somewhat surprisingly, we can do better than the sample kurtosis. In particular, \cite{auddy2023large} introduced an estimate, denoted by $\hat{\scrK}^{\rm AY}_4$, such that
	$$
	\|\hat{\scrK}^{\rm AY}_4-\mathscr{K}_4(\bX)\|=O_p\left(\sqrt{d\over n}\right),
	$$
	under mild regularity conditions. This immediately yields that the eigenvectors, denoted by $\hat{\ba}^{\rm AY}_k$s, of $\hat{\scrK}^{\rm AY}_4$ satisfy
	$$
	\sin\angle(\hat{\ba}_{\pi(k)}^{\rm AY},\ba_k)=O_p\left(\sqrt{d\over n}\right).
	$$
	The estimate $\hat{\ba}^{\rm AY}_k$ is not only consistent for $\ba_k$ whenever $n\gg d$, but also minimax optimal in that no other estimate converges at a faster rate over all possible mixing matrices. Unfortunately, $\hat{\scrK}^{\rm AY}_4$ is not amenable for computation and impractical for high dimensional problems. 
	
	
	\subsubsection{Algorithm Dependent Error Bound} In practice, to estimate $\ba_k$ from the decomposition in Equation~\ref{eq:ica-kurt}, one can use tensor power iteration on the sample kurtosis tensor, starting from some nontrivial initial estimators $\hat{\ba}_k^{[0]}$. This idea was first introduced as cumulant based FastICA by \cite{hyvarinen1999fast} and has since been studied by numerous works in the literature, including \cite{hyvarinen2000independent,hyvarinen2002independent,ollila2009deflation} in the fixed dimension case, and \cite{anandkumar2014tensor,anand2014sample,belkin2013blind,belkin2018eigenvectors,auddy2023large} in the growing dimension case. Given an initialization $\hat{\ba}_k^{[0]}$, the FastICA iteration takes the form
	\begin{equation*}
		\hat{\bw}
		:=
		\hat{\ba}_k^{[t]}
		-
		\dfrac{1}{3n}
		\sum_{i=1}^n\langle\bX_i,\hat{\ba}_k^{[t]}\rangle^3\bX_i
		\,\,\text{;}\,\,
		\hat{\ba}_k^{[t+1]}
		:=\hat{\bw}/\|\hat{\bw}\|.
	\end{equation*}
	Assuming that the initialization is nontrivial in that $\sin\angle(\hat{\ba}^{[0]},\ba_k)<1-\eta$ for some fixed $\eta>0$, \cite{auddy2023large} showed that after $T\gtrsim \log d$ iterations, we have
	$$
	\sin\angle(\hat{\ba}_{\pi(k)}^{[T]},\ba_k)=O_p\left(\sqrt{d\over n}\right).
	$$
	provided $n\gtrsim d^2$. The main challenge is, however, how to obtain a nontrivial initialization.
	
	\subsubsection{Initialization and Computational Gap} A common way to initialize is random slicing. Let $\hat{\scrK}_4$ be an estimator of $\scrK_4(\bX)$. One can initialize via the SVD of: 
	\begin{equation}\label{eq:ica-mat}
		\hat{\mathscr{K}}_4\times_{1,2}\bM 
		=\sum_{k=1}^d\kappa_4(S_k)\langle\ba_k\circ\ba_k,\bM\rangle\ba_k\ba_k^{\top}
		+
		(\hat{\mathscr{K}}_4-\scrK_4)\times_{1,2}\bM. 
	\end{equation}
	Here $\bM\in \RR^{d\times d}$ has iid $N(0,1)$ entries. By sampling $\bM$ repeatedly, one can find $\bM_*$ such that the gap within the singular values $\kappa_4(S_k)\langle\ba_k\circ\ba_k,\bM_*\rangle$ is sufficiently larger than the error $(\hat{\mathscr{K}}_4-\scrK_4)\times_{1,2}\bM_*$. Thus performing an SVD of the matrix in Equation~\ref{eq:ica-mat} provides unit vectors $\hat{\ba}_k^{[0]}$ satisfying $\sin\angle(\hat{\ba}_k^{[0]},\ba_k)<1-\eta$ for some fixed $\eta>0$, after which one can apply the FastICA iterations. Once $\ba_k$ is estimated by $\hat{\ba}_k$, we enforce the slicing matrices $\bM$ to be perpendicular to the already estimated directions, i.e., $\langle \hat{\ba}_k\circ\hat{\ba}_k,\bM\rangle=0$, so that the singular values $\kappa_4(S_k)\langle\ba_k\circ\ba_k,\bM\rangle$ are necessarily small in the directions close to $\ba_k$. 
	
	\cite{anand2014sample} choose $\hat{\scrK}_4=\hat{\scrK}^{\rm sample}_4$ in Equation~\ref{eq:ica-mat}, and show that one has nontrivial initializations whenever $n\gg d^3$.  \cite{auddy2023large} improved this sample complexity through a refined kurtosis estimator that yields a nontrivial estimator whenever $n\gg d^2$. Through low degree polynomial algorithms, \cite{auddy2023large} also provide evidence to show that the lower bound of $n \gg d^2$ is unavoidable in order to obtain a computationally feasible ICA procedure in this setting. Similar to the other tensor related problems, this presents a gap between the sample complexity for information theoretic feasibility, $n\gg d$, and computational tractability, $n\gg d^2$.

	\subsection{Mixture Models}
	
	Tensor-based methods are particularly useful in learning high dimensional mixture models. The classical Gaussian mixture model (GMM) assumes a collection of observations is independently sampled from a common distribution, denoted by ${\rm GMM}\big(\bp, \{\bmu_j\}_{j\in[m]}\big):=\sum_{j=1}^m p_j\cdot N(\bmu_j, \II_d)$, which is a mixture of isotropic $d$-dimensional normal distributions with mean vectors $\bmu_j$s. Here, $p_j$s are the mixing probabilities satisfying ${\bf 1}^{\top}\bp=1$. \cite{hsu2013learning} considered the adjusted second and third order moments, which admit the following decomposition: 
	\begin{align*}
		\bM_2:&=\EE(\bX\circ \bX)-\II_d=\sum_{j=1}^m p_j (\bmu_j\circ \bmu_j)\nonumber\\
		\scrM_3:&=\EE(\bX\circ \bX\circ \bX)-\sum_{i=1}^d \EE\big(\bX\circ \be_i\circ \be_i+\be_i\circ \bX\circ \be_i+\be_i\circ \be_i\circ \bX\big)\\
		&=\sum_{j=1}^m p_j(\bmu_j\circ \bmu_j\circ \bmu_j)\nonumber,
	\end{align*}
	Tensor methods enjoy special advantages in learning GMM. For instance, the eigenvectors of $\scrM_3$ can be uniquely defined even if $p_j=p_k$ for some pair $j\neq k$.  This advantage is unseen in its matrix counterpart $\bM_2$.
	
	\subsubsection{Error Bounds} If $\bmu_1,\cdots,\bmu_m\in\RR^d$ are orthogonal to each other, $\scrM_3$ becomes an ODECO tensor. As in the ICA case, we can derive an estimate for them from an estimate of $\scrM_3$. In particular, let $\hat{\bmu}_k$ be the eigenvectors of an ODECO approximation of the sample moment, denoted by $\hat{\scrM}_3^{\rm sample}$. The perturbation bounds from Equation~\ref{eq:odec-pert2} immediately implies that
	$$
	\sin\angle\left(
	\hat{\bmu}_{\pi(k)},\bmu_k
	\right)
	\le
	\dfrac{\|\hat{\scrM}_3^{\rm sample}-\scrM_3\|}{p_k}
	\quad
	\text{ for all }
	k\in[m]
	,
	$$
	for some permutation $\pi:[m]\to [m]$. Furthermore, it converges to zero in probability if $n\gg d^{3/2}/(\min_{j}p_j)^2$ under suitable moment conditions. This sample complexity can be further improved using the same technique as \cite{auddy2023large}. Specifically, one can construct a moment estimate in the same fashion as they did for estimating the kurtosis tensor, yielding an estimate of $\bmu_k$s satisfying:
	$$
	\sin\angle\left(
	\hat{\bmu}_{\pi(k)},\bmu_k
	\right)
	\lesssim
	\dfrac{1}{p_k}\sqrt{d\over n}
	\quad
	\text{ for all }
	k\in[m]
	$$
	with high probability, attaining a rate which is also minimax optimal. However, as in the ICA case, this estimating procedure is not computationally tractable.
	
	Computationally tractable approaches based on power iteration have also been investigated by \cite{anandkumar2014tensor,anand2014sample} among others. Let $\big\{(\hat p_k, \hat\bmu_k)\big\}_{k\in[m]}$ be the output of tensor power iteration algorithm applied to $\hat\scrM_3^{\rm sample}$ with sequential deflation and a large number of initializations via random slicing. More specifically, by repeatedly sampling $\bv\sim N(\mathbf{0},\II_d)$, one can find $\bv_*$ such that the top singular vector $\hat{\bmu}^{[0]}$ of $\hat{\scrM}_3\times_3\bv_*$ satisfies $\sin\angle(\hat{\bmu}^{[0]},\bmu_k)<1-\eta$ for some fixed $\eta>0$ and some $k\in[m]$, thus acting as a nontrivial initialization for tensor power iteration. It was shown by \cite{anandkumar2014tensor} that, with high probability, there exists a permutation $\pi:[m]\mapsto [m]$ such that 
	$$
	\sin\angle\left(
	\hat{\bmu}_{\pi(k)},\bmu_k
	\right)
	\lesssim \frac{1}{p_k}
	\sqrt{\dfrac{d}{n}}
	\quad {\rm and}\quad 
	\big|\hat p_{\pi(k)}-p_{k} \big|\lesssim\sqrt{\dfrac{d}{n}} 
	\quad
	\text{ for all }
	k\in [m]
	$$
	provided $n\gtrsim d^2/(\min_k p_k)^2$. Once again, we can see the gap in sample complexity.
	


	\subsubsection{Low rank mixture model}
	Matrix-valued observations routinely arise in diverse applications. A low-rank mixture model (LrMM) was proposed by \cite{lyu2023optimal-aos} to investigate the computational and statistical limits in estimating an underlying low-rank structure given a collection of matrix-valued observations. Under the general LrMM, the matrix-valued observations $\bX_1,\cdots,\bX_n\in\RR^{d_1\times d_2}$ are independently sampled from a common distribution, denoted by $\sum_{j=1}^m p_j\cdot N(\bM_j, \II_{d_1}\circ \II_{d_2}) $, which is a mixture of matrix normal distributions with isotropic covariances and low-rank population center matrices $\bM_j$s. By concatenating these matrices into a tensor $\scrX\in\RR^{d_1\times d_2\times n}$, the expectation of $\scrX$,  conditioned on the latent labels,  admits a low-rank Tucker decomposition. 
	Variants of LrMM have appeared in \cite{chen2021learning} for studying low-rank mixed regression, and in \cite{mai2022doubly, sun2019dynamic} for tensor clustering. \cite{lyu2023optimal-aos} showed that a polynomial-time algorithm combining a tensor-based initialization method and a modified second order moment method can achieve a statistically optimal rate in estimating the low-rank population center matrix under the symmetric two-component LrMM. A low-rank Lloyd's algorithm was proposed by \cite{lyu2022lloyd} with a tensor-based initialization method, which achieves a minimax optimal clustering error rate characterized by the separation strength between the population centers. The authors also provided convincing evidence to support the existence of  a statistical-to-computational gap in LrMM for both estimation and clustering. 
	
	\subsection{Tensor Completion}
	Tensor data often have missing values. Tensor completion \citep{liu2012tensor} refers to the problem of recovering a low-rank tensor by observing only a small subset of its entries. Let $\scrT\in\RR^{d\times d\times d}$ have multilinear ranks $(r,r,r)$, and denote $\Omega\subset [d]\times [d]\times [d]$ the locations where the entries of $\scrT$ are observed.   It is commonly assumed that $\mathscr{T}$ is incoherent to ensure the well-posedness of the problem. For notational simplicity, we treat the incoherence parameter of $\scrT$ as a constant, implying that $d^{3/2}\|\scrT\|_{\max}/\|\scrT\|_{\rm F}=O(1)$.  Denote by $n:=|\Omega|$ the sample size. 
	
	\subsubsection{Convex Approaches} In the seminal work by \cite{yuan2016tensor}, a convex program was introduced for exact tensor completion by minimizing the tensor nuclear norm:
	\begin{align}\label{eq:YZ}
		\hat \scrT^{\rm YZ}:=\arg\min\ \|\scrY\|_{\ast},\quad {\rm s.t.}\quad \calP_{\Omega}(\scrY)=\calP_{\Omega}(\scrT),
	\end{align}
	where $\calP_{\Omega}(\scrT)$ is the operator which zeros out the entries of $\scrT$ except those in $\Omega$, and the tensor nuclear norm is defined by $\|\scrY\|_{\ast}:=\sup_{\|\scrM\|\leq 1} \langle \scrY, \scrM\rangle$. Suppose that $\Omega$ is sampled uniformly from $[d]\times [d]\times [d]$ with a cardinality $n\geq C_1\big(r^{1/2}d^{3/2}+r^2d\big)\log^3d$ for a large absolute constant $C_1>0$. \cite{yuan2016tensor} shows that the solution to Equation \ref{eq:YZ} can exactly recover the underlying tensor $\scrT$ with high probability. 
	
	Later, \cite{yuan2017incoherent} introduced the tensor incoherent norms, showing that minimizing the incoherent nuclear norm can, w.h.p., exactly recover an order-$p$ tensor if the sample size $|\Omega|\gtrsim r^{(p-1)/2}d^{3/2}\log^2d$. \cite{ghadermarzy2019near} studied a tensor max-norm and further improved the sample complexity, which nearly matches the information-theoretic lower bound. However, these convex norms are, in general, computationally NP-hard. \cite{barak2016noisy} provided evidence that no polynomial-time algorithms can consistently recover a rank-one $d\times d\times d$ tensor if the sample size $|\Omega|=o(d^{3/2})$.  
	Another line of convex methods is based on semi-definite program (SDP). The sum-of-squares (SOS) method reformulates tensor completion as a large-scale SDP, whose size depends on the degree of relaxation. \cite{potechin2017exact} shows that SOS can exactly recover a $d\times d\times d$ tensor with CP-rank $r$ w.h.p. if the sample size $|\Omega|\gtrsim rd^{3/2}{\rm Polylog}(d)$. Although the SOS method is polynomial-time computable in theory, in practice, its computation is prohibitively intensive due to the need to solve very large-scale SDPs.
	
	\subsubsection{Nonconvex Approaches} Nonconvex methods for tensor completion usually admit much fast computation and consume less memory and storage. In a nutshell,  non-convex methods for completing a $d\times d\times d$ tensor with multilinear ranks $(r,r,r)$ aim to solve
	\begin{align}\label{eq:TC-ncvx}
		\min\ \Big\|\calP_{\Omega}\big([\scrC; \bU_1, \bU_2, \bU_3]-\scrT\big) \Big\|_{\rm F}\quad {\rm s.t.}\quad \scrC\in\RR^{r\times r\times r};  \bU_{1}, \bU_2, \bU_3\in\RR^{d\times r}. 
	\end{align}
	Alternating minimization and gradient descent algorithms are commonly used for finding a locally optimal solution to Equation~\ref{eq:TC-ncvx}.  The success of these algorithms crucially depend on the availability of a warm initialization. Finding a warm initialization is often the most challenging step in tensor optimization problem. Suppose that $\Omega$ is sampled uniformly with replacement. The data essentially consists of a random sample $\{(\scrX_i, Y_i): i\in[n]\}$ of i.i.d. observations satisfying $Y_i=\langle \scrX_i, \scrT\rangle$, where $\scrX_i$ is sampled uniformly from the orthonormal basis $\calX:=\big\{\be_i\circ \be_j\circ \be_k: i, j, k\in[d]\big\}$. \cite{xia2019polynomial} proposed a second-order moment method, which applies spectral initialization to the U-statistic:
	\begin{align}\label{eq:U-stat}
		\frac{d^6}{n(n-1)}\sum_{1\leq i\neq j\leq n}Y_iY_j \calM_k(\scrX_i)\calM_k^{\top}(\scrX_j)\quad \textrm{ for all }k=1,2,3,
	\end{align}
	which is an unbiased estimator of $\calM_k(\scrT)\calM_k^{\top}(\scrT)$. \cite{xia2019polynomial} showed that the objective function in Equation~\ref{eq:TC-ncvx} behaves like a parabola around the oracle and studied the Grassmannian gradient descent algorithm for minimizing the objective function with spectral initialization obtained from Equation~\ref{eq:U-stat}. Let $\hat\scrT^{\rm XY}$ denote the estimator delivered by their algorithm.  \cite{xia2019polynomial} shows that if $n\geq C_1\big(r^{7/2}d^{3/2}\log^{7/2}d+r^7d\log^{6}d\big)$, then $\hat\scrT^{\rm XY}=\scrT$ with high probability.

	\cite{jain2014provable} showed that, once well initialized, the alternating minimization algorithm can exactly recover $\scrT$ with CP-rank $r$ if $n\gtrsim r^5d^{3/2}\log^4 d$. A spectral method was studied by \cite{montanari2018spectral} for both undercomplete ($r\leq d$) and overcomplete ($r\gg d$) regimes, which cannot exactly recover the underlying tensor.

	\subsubsection{Noisy tensor completion} Consider the noisy case where $Y=\langle \scrT, \scrX\rangle+\eps$ has a sub-Gaussian noise $\eps$ with a proxy-variance $\sigma^2$. Noisy tensor completion was first studied by \cite{barak2016noisy} using the sum-of-squares method, but their estimator is statistically sub-optimal. 
	The sample mean $\hat\scrT^{\rm sample}:=d^3n^{-1}\sum_{i=1}^n Y_i\scrX_i$ is an unbiased estimator of $\scrT$. The concentration inequality of $\hat\scrT^{\rm sample}$ was investigated by \cite{yuan2016tensor} and generalized by \cite{xia2021statistically}. For instance, Theorem 1 of \cite{xia2021statistically} shows that $\hat\scrT^{\rm sample}$ concentrates at $\scrT$ in the following manner:
	$$
	\big\|\hat \scrT^{\rm sample}-\scrT \big\|\leq C\alpha\big(\|\scrT\|_{\max}+\sigma\big)\cdot \max\bigg\{\bigg(\frac{d^4}{n}\bigg)^{1/2},\ \frac{d^{3}}{n} \bigg\}\log^5d,
	$$
	which holds with probability at least $1-d^{-\alpha}$ for any $\alpha\geq 1$ and $C>0$ is an absolute constant. 	Note that concentration inequality holds for any given tensor $\scrT$, whether it is full rank or not.  It is worth noting that \cite{xia2021effective} studied a concentration inequality for the case when entries are sampled non-uniformly and without replacement.
	
	While the concentration inequality for the sample mean estimator $\hat \scrT^{\rm sample}$ is sharp in the tensor spectral norm, it is typically full rank, and the error rate it achieves is sub-optimal in the Frobenius norm.   
	The regularized-HOOI algorithm was investigated by \cite{xia2021statistically} showing that the resulting estimator, denoted by $\hat\scrT^{\rm XYZ}$, achieves the minimax optimal rate
	$$
	d^{-3/2}\big\|\hat\scrT^{{\rm XYZ}}-\scrT\big\|_{\rm F}\leq C\big(\|\scrT\|_{\max}+\sigma\big)\cdot \bigg(\frac{(rd+r^3)\log d}{n}\bigg)^{1/2},
	$$
	which holds with high probability if the sample size $n\gtrsim (r^2d+rd^{3/2}){\rm Polylog}(d)$ and the signal strength satisfies
	$\lambda_{\min}/\sigma\gtrsim d^{9/4}n^{-1/2}{\rm Polylog}(d)$. The rate is minimax optimal with respect to the sample size and the degrees of freedom, but is not proportional to noise standard deviation. If $\|\scrT\|_{\max}\gg \sigma$, the rate becomes sub-optimal.
	
	Statistically optimal estimators have also been investigated by \cite{cai2019nonconvex} and \cite{cai2022provable} for noisy tensor completion under the CP and tensor-train model, respectively. A vanilla gradient descent algorithm was studied in \cite{cai2019nonconvex} and \cite{cai2020uncertainty} for completing a CP-format tensor from noisy observations. By a sophisticated leave-one-out analysis, their estimator is shown to be minimax optimal not only in Frobenius norm but also in max-norm.

	\subsection{Tensor Regression}
	
	Tensor regression is a versatile tool with applications in brain imaging , facial image analysis, spatiotemporal learning, multitask learning, chemometrics and myriad other fields. See, e.g., \cite{zhou2013tensor,yu2016learning}.
	
	Let us consider scalar responses depending linearly on tensor covariates: 
	\begin{equation}\label{eq:tens-reg}
		Y_i=\langle\scrX_i,\scrT\rangle+\eps_i,
		\quad
		i=1,\dots,n,
	\end{equation}
	where $Y_i$ is the real valued response, $\eps_i$ is the observation noise for $i=1,\dots,n$. Further $\scrX_i$ and $\scrT$ are the order-3 $d\times d\times d$ dimensional tensors of design and the coefficients respectively, with $\langle\scrX_i,\scrT\rangle$ denoting their element-wise inner product. Equation~\ref{eq:tens-reg} suggests a least squares estimator:
	\begin{equation}\label{eq:tens-reg-ls}
		\hat{\scrT}
		:=
		\underset{\scrA\in \RR^{d\times d\times d}}{\arg\min}
		\sum_{i=1}^n
		(Y_i-\langle\scrX_i,\scrA\rangle)^2.
	\end{equation}
	Stringing $\scrA$ into a $d^3$ dimensional vector and using standard regression methods is statistically inefficient. Instead we now describe some methods which benefit from the multiway nature of $\scrA$ and take advantage of its low dimensional structures. 
	
	\subsubsection{Convex approaches} A common approach of solving the constrained least squares problem in Equation~\ref{eq:tens-reg-ls} is to add a convex penalty:
	\begin{equation}\label{eq:tens-reg-pen}
		\hat{\scrT}
		:=
		\underset{\scrA\in \RR^{d\times d\times d}}{\arg\min}
		\left\{
		\dfrac{1}{2n}
		\sum_{i=1}^n(Y_i-\langle\scrX_i,\scrA\rangle)^2
		+\tau\calR(\scrA)
		\right\}.
	\end{equation}

	The penalty $\calR(\cdot)$ is chosen based on the intended structural constraints. For elementwise sparsity, one uses $\calR(\scrT)=\|{\rm vec}(\scrT)\|_1=\sum_{i,j,k}|\scrT_{ijk}|$. On the other hand for the sparsity of fibres, we use a group-based regularizer $\calR(\scrT)=\sum_{i,j}\|\scrT_{ij.}\|_2$, akin to the group Lasso.	It is also common to make low rank assumptions, e.g., via a low CP-rank $r$: which can be imposed by a tensor nuclear norm penalty $\calR(\scrT)=\|\scrT\|_*=\max_{\mathscr{S}:\|\mathscr{S}\|=1}\langle\scrT,\scrS\rangle$. Low multirank assumptions can be placed similarly by penalizing through the sum of the nuclear norms of the matricizations $\sum_q\|\calM_q(\scrT)\|_*$.
	
	
	\cite{raskutti2019convex} provided a general framework for analyzing convex regularized tensor regression estimators of the form Equation~\ref{eq:tens-reg-pen} with \emph{weakly decomposable penalties} $\calR(\cdot)$. Within this framework, they defined two key quantities:
	$$
	s(\calA)=\sup_{\scrA\in\calA/\{0\}}\dfrac{\calR(\scrA)^2}{\|\scrA\|_{\rm F}^2}
	\text{ and }
	w_G(\calS)=\EE\left(\sup_{\scrA\in\calS}\langle\scrA,\scrG\rangle\right)
	$$
	where $\scrG\in\RR^{d\times d\times d}$ is a tensor with independent $N(0,1)$ entries. Here, roughly speaking, $\calA$ denotes the low dimensional subspace of tensors in which we constrain $\scrT$ to be in. The first quantity $s(\calA)$ then measures the intrinsic dimensionality of $\calA$. The second quantity, called Gaussian width, is a notion of the complexity of $\calS$. When evaluated on $\BB_{\calR}(1):=\{\scrA:\calR(\scrA)\le 1\}$, this second quantity determines how large the penalty $\calR$ is with respect to the Frobenius norm. With these definitions, \cite{raskutti2019convex} proved that the solution $\hat{\scrT}$ of Equation~\ref{eq:tens-reg-pen}, with high probability satisfies
	$$
	\|\hat{\scrT}-\scrT\|^2_{\rm F}\lesssim s(\calA)\tau^2
	$$
	provided $\tau\gtrsim \sigma w_G[\BB_{\calR}(1)]/\sqrt{n}$, under Gaussian designs and errors $\eps_j\sim N(0,\sigma^2)$. 
	
	When specialized to elementwise sparsity of the form $\sum_{ijk}\mathbbm{1}(\scrT_{ijk}\neq 0)\le s$, the general bound above boils down to $\|\hat{\scrT}-\scrT|_{\rm F}^2\lesssim \sigma^2s\log(d)/n$, provided $\tau\asymp\sigma\sqrt{\log(d)/n}$. When imposing a low rank structure, say through the multirank condition $\max\{{\rm rank}(\calM_q(\scrT)):q=1,2,3\}\le r_{\max}$ one can choose the penalty $\calR(\scrT)=\sum_q\|\calM_q(\scrT)\|_*$. Then the solution $\hat{\scrT}$ of Equation~\ref{eq:tens-reg-pen} satisfies $\|\hat{\scrT}-\scrT\|_{\rm F}^2\lesssim \sigma^2r_{\max}d^2/n$, by taking $\tau\asymp \sigma \sqrt{d^2/n}$.

	
	\subsubsection{Nonconvex approaches} While convex relaxations to the nonconvex tensor regression problem provide polynomial time algorithms, they are often slow in practice. Moreover, a convex relaxation might involve quantities computing which are NP-hard. The nuclear norm $\|\scrT\|_*=\sup\{\langle\scrA,\scrT\rangle:\|\scrA\|\le 1\}$, used to impose CP low rank-ness in tensor regression, is one such example. In stead, nonconvex projected gradient descent algorithms for tensor regression are often beneficial for both computational and statistical purposes. 
	
	\cite{chen2019non} consider a generalized linear model for a scalar response $Y$ and a covariate tensor $\scrX\in\RR^{d_1\times d_2\times d_3}$, where the conditional distribution of $Y$ on $\scrX$ is given by 
	$
	p(Y|\scrX,\scrT)=h(Y)\exp\{Y\langle\scrX,\scrT\rangle-a(\langle\scrX,\scrT\rangle)\}
	$
	where $a(\cdot)$ is strictly convex and differentiable log-partition function. This results in the negative log-likelihood objective function
	$$
	\calL(\scrA)=\dfrac{1}{n}\sum_{i=1}^n\left[
	a(\langle\scrX_{i},\scrA\rangle)-Y_i\langle\scrX_i,\scrA\rangle
	-\log h(Y_i)
	\right].
	$$
	In this subsection, we consider projected gradient descent (PGD) to optimize the above objective function. The PGD framework, introduced by \cite{jain2014iterative}, can be applied to a general tensor space  $\calF$ as follows: 
	$$
	\hat{\scrT}^{[t+1]}\longleftarrow \textsf{Proj}_{\calF}\Big(\hat \scrT^{[t]}-\eta\nabla \calL\big(\hat \scrT^{[t]}\big)\Big).
	$$
	In practice $\calF$ can be described based on low rank or sparsity assumptions. In \cite{chen2019non} the authors derive theoretical guarantees of the PGD estimator for general loss functions satisfying restricted strong convexity and for collections of subspaces that allow for contractive projections. Due to the nonconvexity, the performance of PGD depends on the initialization. This is reflected in the local Gaussian width $w_G[\calF_0\cap \BB_{\rm F}(1)]$. Here $\calF_0$ denotes the subspace where the PGD iteration is initialized, and $\BB_{\rm F}(1)$ is the unit ball of tensors in Frobenius norm. With these conditions, \cite{chen2019non} show that:
	$$
	\|\hat{\scrT}^{[k]}-\scrT\|_{\rm F}\le c w_G[\calF_0\cap\BB_{\rm F}(1)]/\sqrt{n}+\delta
	$$ 
	for some constant $c>0$, provided PGD is run for $k\gtrsim \log(\|\scrT\|_{\rm F}/\delta)$ iterations. Thus if we have $n\gtrsim (w_G[\calF_0\cap\BB_{\rm F}(1)])^2$ samples, the PGD estimator converges linearly to the true covariate tensor $\scrT$. 
	
	Specializing to the linear model and the sparsity/low rank examples considered earlier, one can compare the nonconvex and convex approaches. For illustrative purposes, we only consider the low Tucker rank regression model, i.e., $\scrT$ has multirank at most $(r_1,r_2,r_3)$ with $r_{\min}:=\min\{r_1,r_2,r_3\}$. In this setting, \cite{chen2019non} used an approximate projection operator defined through a low rank projection on the matricizations. With this estimator, they showed that with $n\gtrsim r_{\min}d^2$ samples, the PGD estimator $\hat{\scrT}$ satisfies $\|\hat{\scrT}-\scrT\|_{\rm F}^2\lesssim \sigma^2r_{\min}d^2/n$ after $\log(n)$ many iterations. The PGD estimator improves on the convex relaxed estimator of \cite{raskutti2019convex} by selecting the optimal mode for matricization, whereas convex regularization takes an average of the three matricizations, and is hence sub-optimal. 
	
	Following the above works, in \cite{zhang2020islet} the authors addressed the issue of scalability of tensor regression to massive dimensions. Their method is based on importance sketching (ISLET) where, depending on low rankness or sparsity requirements, HOOI or STAT-SVD (see \cite{zhang2019optimal}) is used to determine the importance sketching directions. Then tensor regression is applied on the dimension reduced tensors obtained from importance sketching. Due to the dimension reduction, ISLET has immense computational as well as theoretical advantage over a regression method applied on the entire tensor. 

	\subsection{Higher Order Networks}
	Higher order networks model the joint interactions of multiple entities, such as hypergraph networks \citep{benson2016higher}, multi-layer networks \citep{kivela2014multilayer}, and dynamic networks \citep{mucha2010community}. These higher-order networks can be conveniently represented as higher-order tensors. Consider a $3$-uniform hypergraph on $d$ vertices, where  each hyperedge consists of exactly three vertices. It is equivalently representable by the adjacency tensor $\scrA\in\{0,1\}^{d\times d\times d}$. Specifically, the entry $A_{i_1i_2i_3}=1$ if and only if the hyperedge $(i_1,i_2,i_3)$ is present. If the hyperedges are undirected, tthen he adjacency tensor is symmetric.  
	
	Clustering nodes in an undirected hypergraph network, often referred to as community detection, is of particular interest in many applications. Suppose there exist $K$ groups among $d$ vertices, characterized by a binary matrix $\bZ\in\{0,1\}^{d\times K}$, where each row has exactly one entry being $1$ and the others are $0$s. Under the hypergraph stochastic block model (hSBM, \cite{ghoshdastidar2014consistency}), the hyperedge connectivity is determined solely by the node community memberships, summarized by a probability tensor $\scrP\in[0,1]^{K\times K\times K}$. The expected  adjacency tensor $\scrA$ admits a Tucker decomposition as 
	\begin{align*}
		\EE[\scrA]=[\scrC; \bZ, \bZ, \bZ]. 
	\end{align*}
	Since the singular vectors of $\EE[\scrA]$ provide community membership information, one natural idea is to estimate the population singular vectors by finding a low-rank approximation of $\scrA$.  If the oracle cluster sizes are balanced, the singular vectors of $\EE[\scrA]$ are incoherent. Motived by this, \cite{ke2019community} proposed to find a regularized tensor decomposition of $\scrA$ based on the regularized-HOOI algorithm. They derived a concentration inequality of $\scrA$ using the incoherent operator norm:
	\begin{align*}
		\|\scrA-\EE[\scrA]\big\|_{o,\delta}:=\sup_{\substack{\ba,\bb,\bc\in\SS^{d-1} \\ \|\bc\|_{\max}\leq \delta}} \big<\scrA-\EE[\scrA], \ba\circ \bb\circ \bc\big>
	\end{align*}
	for any $\delta\in(0,1]$. The conventional tensor operator norm corresponeds to the special case $\|\cdot\|_{o,1}$. 
	If the cluster sizes are balanced, one can focus on small values of $\delta$s, such as $O\big((K/d)^{1/2}\big)$.  This provides technical convenience for dealing with extremely sparse hypergraph networks.  
	
	Suppose that the oracle clusters have balanced sizes and $\scrP/q_e$ is full rank with a bounded condition number. Denote $q_e:=\|\scrP\|_{\max}$, which characterizes the overall network sparsity level . Indeed, $d^2q_e$ can be viewed as the average node degree. \cite{ke2019community} shows that if $d^2q_e\gg K^2\log^2d$, then the regularized regularized-HOOI and spectral clustering algorithm can consistently recover all the communities with probability ending to one as $d\to\infty$. It means that their method is effective as long as 
	the total number of hyperedges exceeds $d\log^3d$, assuming $K=O(1)$. Therefore, the regularized tensor decomposition method can handle extremely sparse hypergraph networks. It is worth noting that the method studied in \cite{ke2019community} can also handle node heterogeneity. 
	
	Regularized tensor decomposition has been investigated by \cite{jing2021community} for community detection in multi-layer networks and by \cite{lyu2023latent} for general higher-order network analysis. In particular, a sharp concentration inequality for the adjacency tensor with respect to incoherent operator norm was developed by \cite{jing2021community}. They show that the regularized-HOOI and spectral clustering can consistently recover both the global and local community memberships for multi-layer networks. Their theories allow the network sparsity to decrease when more layers are observed. The methods are still effective even when some layers are disconnected graphs. More recently, \cite{lyu2023optimal} studied a two-step algorithm for clustering layers in a multi-layer network and achieved the minimax optimal clustering error rate under a wide range of network sparsity conditions.
	
	The tensor block model (TBM) can be viewed as a generalization of hSBM. A rank-one sparse TBM was proposed by \cite{xia2019sup}, assuming that the observed tensor $\scrA$ satisfies $\scrA=\lambda {\bf 1}_{C_1}\circ {\bf 1}_{C_2}\circ {\bf 1}_{C_3}+\scrE$, where $C_j\subset [d]$, and the noise tensor $\scrE$ has i.i.d. centered sub-Gaussian entries. Here, ${\bf 1}_C$ denotes the binary vector with only the entries in $C$ being equal to one.  More general TBMs have been investigated by \cite{wang2019multiway}, \cite{han2022exact} and \cite{zhou2023heteroskedastic} for cluster analysis. In particular, \cite{han2022exact} proposed a higher-order spectral clustering algorithm for finding clusters in each dimension, where a statistical-to-computational gap was discovered. Their algorithm can exactly recover all the clusters under a suitable separation condition.

	\subsection{Tensor Time Series}
	
	Tensor time series have become ubiquitous to date. For example, the monthly international trade flow of different categories of commodities between countries produces a time series of tensor-valued observations \citep{cai2022generalized,chen2022factor}. As elaborated by \cite{rogers2013multilinear}, the challenge in tensor time series analysis is to study the dynamics while preserving the underlying tensor structures.
	
	Tensor factor model \citep{chen2022factor} assumes the tensor-valued observations share common loadings and the dynamic is characterized by smaller-sized factor tensors, which is generalized from matrix factor model \citep{chen2019constrained,wang2019factor}. Let $\scrX_t\in\RR^{d_1\times d_2\times d_3}, t\in[T]$ be a time series of tensors.  Under tensor factor model, it is assumed that 
	$$
	\scrX_t=\big[\scrF_t; \bU_1, \bU_2, \bU_3\big]+\scrE_t,
	$$
	where $\scrF_t\in\RR^{r_1\times r_2\times r_3}, t\in[T]$  are factor tensors. The Tucker decomposition preserves the tensor structure, while the factor tensor time series save the underlying dynamic. It is assumed that the process $\scrF_t$ is weakly stationary and the lagged cross-products have fixed expectation in that $(T-h)^{-1}\sum_{t=h+1}^{T} \scrF_{t-h}\circ \scrF_t$ converges to $\EE [\scrF_{T-h}\circ \scrF_{T}]$ in probability. For simplicity, the tensor $\scrE_t$s are assumed to be white noise and mutually independent across time. Estimating procedures under tensor factor model can be technically involved. \cite{chen2022factor} proposed two estimators, named TOPUP (time series outer-product unfolding procedure) and TIPUP (time series inner-product unfolding procedure). For any $k=1,2,3$, TOPUP and TIPUP calculate the mode-$k$ lag-$h$ auto-correlations using outer and inner products as follows:
	$$
	\bV_{k,h}^{(O)}:=\frac{1}{T-h}\sum_{t=h+1}^{T} \calM_k(\scrX_{t-h})\circ \calM_k(\scrX_t)\,\, {\rm and}\,\, \bV_{k,h}^{(I)}:=\frac{1}{T-h}\sum_{t=h+1}^{T} \calM_k(\scrX_{t-h})\calM_k^{\top}(\scrX_{t}),
	$$ 
	respectively. Conditioned on $\{\scrF_t: t\in[T]\}$, the left singular subspaces of $\EE \bV_{k,h}^{(O)}$ and $\EE \bV_{k,h}^{(I)}$ are both directly related to the loading matrix $\bU_k$. TOPUP and TIPUP then exploits multiple auto-correlations differently where the former  and latter ones take the left singular vectors of the  assembles
	$$
	\big(\calM_1(\bV_{k,h}^{(O)}), h=1,\cdots, h_0\big)\quad {\rm and}\quad \big(\calM_1(\bV_{k,h}^{(I)}), h=1,\cdots, h_0\big),
	$$
	respectively. Finite sample properties of TOPUP and TIPUP were established by \cite{chen2022factor}. A CP-format tensor factor model was proposed by \cite{han2021cp}. The authors introduced a novel initialization procedure based on composite PCA and proposed to refine the estimates by iterative simultaneous orthogonalization. The computational and statistical performances of their algorithm were investigated assuming that the factor time series is stationary and strongly $\alpha$-mixing.

	Autoregressive (AR) models are useful for predicting future observations in time series. Matrix and tensor time series AR models have been proposed by \cite{chen2021autoregressive} and \cite{wang2024high}. For instance, under a matrix time series AR$(q)$ model, the observations $\bX_t\in\RR^{d_1\times d_2}, t\in[T]$ are assumed to satisfy
	$$
	\bX_t=\bA_1\bX_{t-1}\bB_1^{\top}+\cdots+\bA_q\bX_{t-q}\bB_q^{\top}+\bE_t,\quad t=q+1,\cdots, T,
	$$
	where $q\geq 1$, $\bA_k$s and $\bB_k$s are unknown $d_1\times d_1$ and $d_2\times d_2$ matrices, and $\bE_t$s are stochastic noise. Under a tensor AR$(1)$ model, the observations $\scrX_t\in\RR^{d_1\times d_2\times d_3}, t\in[T]$ are assumed to satisfy
	$$
	\scrX_{t}=\scrX_{t-1}\times_1 \bC_1\times_2 \bC_2\times_3 \bC_3+\scrE_t,\quad t=2,\cdots, T,
	$$
	where the unknown matrix $\bC_k$ is of size $d_k\times d_k$. \cite{chen2021autoregressive} and \cite{wang2024high} investigated the least squares method to estimate the unknown matrix parameters and proposed to solve it by resorting to the nearest Kronecker product problem or by alternating minimization. Both estimating procedures are non-convex, and only locally optimal solutions can be found. Asymptotic normality was established, under mild conditions, for the globally optimal solution to the least squares minimization. 
	
	\section{Concluding Remarks}
	Tensors methods are ubiquitous in modern statistical applications. In many instances, the nominal complexity of the problem is much greater than the information content of the data and the main challenge is  to develop inferential methods to achieve both statistical and computational efficiencies. Over the past decade, considerable strides have been taken in the development of innovative statistical techniques, efficient computational algorithms, and fundamental mathematical theory to analyze and exploit information in these types of data. This review samples various representative advancement to assess the statistical and computational performance of the maximum likelihood and other optimization methods, convex regularization,  gradient decent and power iteration, tensor unfolding, and their combinations in a wide range of tensor-related problems. A recurring theme is the intricate interplay between two types of efficiencies -- statistical and computational -- with a gap often existing between the performances of statistically and computationally efficient methods. A deep understanding of this phenomenon has profound significance in statistical theory.

	%
	\bibliographystyle{plainnat}
	\bibliography{review-refs}

\end{document}